\documentclass[12pt,eqno]{article} 

\usepackage{amsmath}
\usepackage{amssymb}

\addtocounter{section}{-0}
\newtheorem{guess}{Theorem}[section]
\newtheorem{proposition}[guess]{Proposition}

\newtheorem{remark}[guess]{Remark}
\newtheorem{lemma}[guess]{Lemma}
\newtheorem{corollary}[guess]{Corollary}

\numberwithin{equation}{section}

\newcommand{\be}{\begin{equation}}
\newcommand{\ee}{\end{equation}}

\newcommand{\jm}{\jmq}
\newcommand{\jmq}{Q^{\alpha , \beta }_N}

\newcommand{\Ninf}{$N\to\infty$}
\newcommand{\lab}{$L(b^N,a^N)$}

\begin{document}

\title{From Toda to KdV}

\author{D. Bambusi\footnote{Dipartimento di Matematica, Universit\`a degli Studi di
Milano, Via Saldini 50, I-20133 Milano}, T. Kappeler\footnote{Institut f\"ur Mathematik,
Universit\"at Z\"urich, Winterthurerstrasse 190, CH-8057 Z\"urich. Supported in part by
the Swiss National Science Foundation.}, T. Paul\footnote{CNRS
and CMLS, \'Ecole Polytechnique, F-91128 Palaiseau  }}

\maketitle

\small

\begin{abstract}
\noindent
For periodic Toda chains with a large number $N$ of particles we
consider states which are $N^{-2}$-close to the equilibrium and
constructed by discretizing arbitrary given $C^2-$functions with mesh
size $N^{-1}.$ Our aim is to describe the spectrum of the Jacobi
matrices $L_N$ appearing in the Lax pair formulation of the dynamics
of these states as $N \to \infty$. To this end we construct two Hill
operators $H_\pm$ -- such operators come up in the Lax pair
formulation of the Korteweg-de Vries equation -- and prove by methods
of semiclassical analysis that the asymptotics as $N \rightarrow
\infty $ of the eigenvalues at the edges of the spectrum of $L_N$ are
of the form $\pm (2-(2N)^{-2} \lambda ^\pm _n + \cdots )$ where
$(\lambda ^\pm _n)_{n \geq 0}$ are the eigenvalues of $H_\pm $.  In
the bulk of the spectrum, the eigenvalues are $o(N^{-2})$-close to the
ones of the equilibrium matrix. As an application we obtain
asymptotics of a similar type of the discriminant, associated to
$L_N$.
\end{abstract}



\section{Introduction}
\label{1. General Introduction}

\noindent It is well known that  the (periodic) Toda lattice 
is an integrable system and by classical {\em heuristic}
arguments, its dynamics are expected to be well described by solutions
of the (periodic) KdV equation in the continuous limit (cf \cite{ZK,C,T}).  However, only
quite recently \cite{SchW, BP}, it has been rigorously proved that in
an appropriate asymptotic regime, small solutions of Toda lattices, or
more generally, of chains of particles with nearest neighbors
interaction, referred to as FPU chains, can be described approximately
in terms of solutions of the KdV equation.  It is important to note
that in order to approximate {\it one} solution of an FPU chain, {\it
  two} solutions of the KdV equation are needed, one corresponding to
a right moving wave, and the other corresponding to a left moving
wave.  Furthermore we recall that these results are proved by
averaging type methods, allowing to control the dynamical variables
for long, but finite intervals of time.

Since the (periodic) Toda lattice is an integrable system, one expects
that the above approximation results can be improved in this case
by computing the asymptotics of quantities such as the frequencies.  
In addition, one would like to understand how the integrable structure of
Toda lattices is related to the corresponding one of the KdV equation 
in the continuous limit. 
In particular, recall that both equations admit a Lax pair formulation. The one 
for Toda lattices involves  Jacobi matrices and the one
for the KdV equation Schr\"odinger operators. 
In the setup of lattices on the entire line, Toda \cite{T} 
showed by a formal computation that the continuous limit of
 Jacobi matrices is given by {\it one} Schr\"odinger operator (cf
 \cite{T}, p.~93), leading to {\it one} solution of the
 KdV equation.  But in view of the rigorous results of \cite{SchW,BP},
 {\it two} solutions of KdV are needed to describe the asymptotics of
 solutions of Toda lattices in the continuous limit.  Hence even on a
 formal level, Toda's result is incomplete, at least for general
 initial data. These limitations of Toda's result are also shared by
 works in the periodic setup such as \cite{Gie1,Gie2} as well as 
by studies of other lattices (cf e.g. \cite{CP}).
Indeed, in \cite{Gie1}, p.~ 587, the author points out that he only considers 
very special initial data of the periodic Toda lattice.

The formal results of Toda et al. and the rigorous results of \cite{SchW, BP}
lead to the problem of how to construct two Schr\"odinger operators
yielding the two KdV solutions needed to describe the asymptotics of Toda lattices 
in the continuous limit without the restrictions on the initial data mentioned above.
In the present paper, we solve this problem in the periodic setup in which case
 these operators are also referred to as Hill operators. It might come as a surprise that 
they are constructed by some methods of semiclassical analysis.
One of the main results we show says that they can be used to approximately describe the limiting 
asymptotics of the spectra of periodic Jacobi matrices. 

We believe that our results and the methods developed for proving them
will be an essential tool for studying all kinds of properties of the
asymptotics of Toda lattices in the continuous limit. Results in this
direction are obtained in \cite{bkp2} by applying what is proved in
the present paper: one of
the results of \cite{bkp2} provides the first two terms in the
asympotics of the Toda frequencies in terms of the KdV frequencies
corresponding to the {\em two} Hill operators, mentioned above.

Finally we would like to discuss the connection of the present
research with the so called FPU problem, which actually is the main
motivation for our research.  We recall that in their celebrated
report \cite{FPU}, Fermi, Pasta, and Ulam studied the dynamics of FPU
chains.  They wanted to confirm numerically that energy sharing among
the different degrees of freedom occurs. Very much to their
surprise, they observed recurrent dynamics instead. It led to the question,
referred to as the FPU problem, whether the observed recurrence phenomenon
persists in the thermodynamic limit.  In case it does, it would
contradict the so called equipartition principle leading to potentially serious
problems for the foundations of classical statistical mechanics.
We emphasize that, notwithstanding the
huge number of computations and the enormous amount of theoretical
work done up till now (see e.g. \cite{BGG}, \cite{BKL}, \cite{LPV},
\cite{PL}), an answer to this problem is still not known. For status
reports on the research of FPU chains see \cite{BI}, \cite{CGG},
\cite{G}.

In the context of the FPU problem, our interest in Toda lattices stems
from the facts that on the one side, Toda lattices are integrable and
hence their continuous limits might be easier to study and at a deeper
level, due to the additional structures present, and that on the other
side, near the equilibrium, FPU chains are well approximated by
Toda chains.  Indeed, it was already pointed
out in \cite{FFM} that close to the equilibrium, (periodic) FPU chains
are typically better approximated by (periodic) Toda chains than by
the linear model.  Subsequently, numerical evidence was found that up
to a long time, periodic solutions of FPU chains with small initial
data are very well approximated by Toda chains.  See the quite recent
work in this direction \cite{BCP} and references therein.


\section{Statement of main results}
\label{1. Introduction}

The Toda lattice, in the setting
of periodic boundary conditions with period $N\geq 2,$ is the
Hamiltonian system with Hamiltonian
 \[
 \mathcal H^{\mbox{}}=\frac 1 2\sum_{n=1}^Np_n^2+\sum_{n=1}^Ne^{q_n-q_{n+1}}.
 \]
 Here $q_n$ denotes the displacement from the equilibrium position of the $n$'th particle, $p_n$ its momentum and $(q_n,p_n)$ is defined for any $n$ in $\mathbb Z$ by requiring that $(q_{i+N},p_{i+N}) =  (q_i,p_i)$ for any $i\in\mathbb Z$.
 When expressed in Flaschka coordinates,
 $
 b_n=-p_n\mbox{ and } a_n=e^{\frac 12(q_n-q_{n+1})}
 $(\cite{F}),
 the Hamiltonian equations of motion associated to $\mathcal H^{\mbox{}}$ take 
a Lax pair formalism description given
by 
\begin{equation}\label{laxx}
\dot L = [B,L]
\end{equation}
 where the $N \times N$ matrices $L = L(b,a)$ and $B = B(a)$ are
of the form
   \[  \begin{pmatrix} 
   b_1 &a_1 &0&\ldots&0&a_N \\ 
   a_1 &b_2 &a_2&\ldots &0&0 \\ 
   0 &a_2 &b_3&\ldots&0&0\\
   \vdots&&&&&\vdots\\
   0&\ldots &&\ldots &b_{N-1}&a_{N-1}\\
      a_N &0&\vdots &&a_{N-1}&b_N 
      \end{pmatrix}
      \mbox{ and } 
      \begin{pmatrix} 
   0 &a_1 &0&\ldots&0&-a_N \\ 
   -a_1 &0 &a_2&\ldots &0&0 \\ 
   0 &-a_2 &0&\ldots&0&0\\
   \vdots&&&&&\vdots\\
   0&\ldots &&\ldots &0&a_{N-1}\\
      a_N &0&\vdots &&-a_{N-1}&0 
      \end{pmatrix}
       \]      
respectively, with $a = (a_n)_{1 \leq n \leq N} \in {\mathbb
  R}^N_{>0}$ and $b = (b_n)_{1 \leq n \leq N} \in {\mathbb
  R}^N$. Notice that the matrix $L(0_N,1_N)$ with $b = 0_N = (0,
... ,0)$ and $a = 1_N = (1, ... , 1)$ is an equilibrium for
\eqref{laxx}. We are interested in the $N \rightarrow \infty $
asymptotics of various spectral quantities of $L(b^N,a^N)$ for
$b^N,a^N$ of the form 
\be\label{Da.300} b^N_n = \varepsilon \beta
\big( \frac{n}{N} \big) \quad \mbox{ and } \quad a^N_n = 1 +
\varepsilon \alpha \big( \frac{n}{N} \big) 
\ee 
where $\varepsilon $ is
a {\it coupling parameter} and $\alpha , \beta$ are functions in
$C_0^2({\mathbb T}, {\mathbb R}),$ i.e. $1-$periodic $C^2-$functions
with $[\alpha ] = [\beta ] = 0$ with $[\alpha ]$ denoting the mean of
$\alpha , [\alpha ] = \int ^1_0 \alpha (x)dx.$ 
Alternatively, one can consider
\[
p^N_n = - \varepsilon \beta \big( \frac{n}{N} \big) \quad \mbox{and} \quad 
q^N_n = - {2N\varepsilon} \xi \big( \frac{n}{N} \big)
\]
where $\xi$ is the element in $C_0^3(\mathbb T)$, satisfying $\xi ' = \alpha$.
Using that
\[
\exp( \frac{q^N_{n} - q^N_{n+1}}{2} ) = 1 + \frac{q^N_{n} - q^N_{n+1}}{2} +O(\epsilon/N) 
= a_n^N + O(\epsilon / N)
\]
one can show that our results stated below hold for either of the two
discretizations.

The limiting equations strongly depend on the choice of the coupling
parameter $\epsilon$. In \cite{BGPU} it is shown that with $\epsilon
\sim 1$, one obtains in the limit as $N \to \infty$ a nonlinear system of
equations of hyperbolic type, which contains as a special case the inviscid Burgers equation.
 In contrast to \cite{BGPU}, we choose $ \varepsilon \equiv \varepsilon _N = (2N)^{-2}.$
It turns out that in this case, the asymptotics of the dynamics is described 
in terms of two solutions of the KdV equation (cf \cite{bkp2}). 
Our aim is to compute the asymptotics of the eigenvalues of
$L(b^N, a^N)$ and of the corresponding discriminant as $N \rightarrow
\infty $.  Let us note that in view of the Lax pair representation,
the spectrum of $L(b^N,a^N)$ is conserved by the Toda flow. To obtain
a set of independent integrals of motion it turns out to be more
convenient (see e.g.  \cite{HK1}) to double the size of \lab\ and to
consider \be\label{jm}Q^{\alpha , \beta }_N \equiv Q(b^N, a^N) = L
\big( (b^N, b^N), (a^N, a^N) \big) , \ee namely
\[
Q^{\alpha , \beta }_N=
 \begin{pmatrix} 
   b^N_1 &a^N_1 &0&\ldots&\ldots&\ldots&\ldots&\ldots&0&a^N_N \\ 
   a^N_1 &b^N_2 &a^N_2&0 &\ldots&\ldots&\ldots&\ldots&\ldots&0 \\ 
   0 &a^N_2 &b^N_3&a^N_3&0&\ldots&\ldots&\ldots&\ldots&0\\
   \vdots&&&&&&&&&\vdots\\
   0&\ldots&0&a^N_{N-1}&b^N_N&a^N_N&0&\ldots&\ldots&0\\
   0&\ldots&\ldots&0&a^N_N&b^N_1&a^N_1&0&\ldots&0\\
   \vdots&&&&&&&&&\vdots\\
   0&\ldots &\ldots&\ldots &\ldots&\ldots&0&a^N_{N-2}&b^N_{N-1}&a^N_{N-1}\\
      a^N_N &0&\ldots &\ldots&\ldots&\ldots&\ldots&0&a^N_{N-1}&b^N_N 
      \end{pmatrix}
\]   
The eigenvalues of $Q^{\alpha , \beta }_N$ when listed in increasing order
and with multiplicities satisfy
\[ \lambda ^N_0 < \lambda ^N_1 \leq \lambda ^N_2 < \cdots < \lambda ^N_{2N-3} \leq
      \lambda ^N_{2N-2} < \lambda ^N_{2N-1} .
\]
By Floquet theory (cf. e.g. \cite{HK1}) one sees that in the case where $N$ is even, $\lambda_0,\lambda_3,\lambda_4,\dots,\lambda_{2N-5},\lambda_{2N-4},\lambda_{2N-1}$ are the $N$ eigenvalues of \lab . For $N$ odd, they are given by
$\lambda_1,\lambda_2,\lambda_5,\lambda_6,\dots,\lambda_{2N-5},\lambda_{2N-4},\lambda_{2N-1}$.
 To describe the asymptotics of $\lambda ^N_n$ at the edges, $n \sim 1$ or $n \sim 2N-1$,
we need to introduce two Hill operators $H_\pm := - \partial ^2_x + q _\pm$ with potentials
   \begin{equation}\label{pot} q_\pm (x) = - 2\alpha (2x) \mp \beta (2x).
   \end{equation}
The discovery of these two operators and of their role in the description of the asymptotics 
as $N \to \infty$ of the spectrum of $Q^{\alpha , \beta }_N$ is one of the main contributions of this paper.
The role played by the operators $H_-$ and $H_+$ in the description of the asymptotics of the left respectively right edge of the spcetrum of $Q^{\alpha , \beta }_N$  will be explained in detail in Section \ref{quasiedges}.
Furthermore we point out that in our subsequent paper \cite{bkp2} we prove that 
the limiting dynamics of the Toda lattices 
$(b^N,a^N)$ as $N \to \infty$ can be described in terms of the solutions of the KdV equation 
corresponding to $q_-$ and $q_+.$ Note that the two potentials $q_-$ and $q_+$ determine $\alpha$
and $\beta$ uniquely and that they are independent from each other. Loosely speaking, in terms of the
asymptotics described in \cite{bkp2}, it means that the parts of the Toda lattices corresponding 
to the spectrum at the two edges do not interact although Toda lattices are nonlinear systems.
The formulas \eqref{pot} for $q_-$ and $q_+$ are an outcome of semiclassical analysis, discussed
in Section \ref{toplitz} -- see also the explanations below after Remark \ref{Remark2.1} .

   Note that $q_\pm$ are periodic functions of period $1/2$. The periodic  eigenvalues $(\lambda ^\pm_n)
_{n \geq 0}$ of $H_\pm $ on $[0,1]$, when listed in increasing order and with multiplicities, are  known to satisfy $\lambda ^\pm _0 < \lambda ^\pm _1 \leq \lambda ^\pm _2 < \cdots $.
It turns out that the asymptotics of the eigenvalues of $\jmq$ exhibit three different regions: the bulk
and the two edges, which shrink to $\{-2\}$ and $\{+2\},$ respectively, as \Ninf. Each of these three parts of the spectrum has its proper asymptotics:  in the bulk, the spectrum is close to the one of the equilibrium matrix by a distance smaller than the distance between the given Jacobi matrices and the equilibrium matrix, whereas in each of the two edges, the first correction is of the same order as this distance and involves the spectrum of one of the two Hill operators $H_{\pm}$.

To define the two edges of the spectrum consider a function 
$F : {\mathbb N} \rightarrow {\mathbb R}_{\geq 1}$
satisfying
\[ 
(F)\qquad \lim _{N \rightarrow \infty} F(N) = \infty; \qquad \mbox{ increasing; } \qquad  F(N) \leq N^\eta \mbox{ with } 0 < \eta < 1/2. \qquad
\]

\begin{guess}
\label{Theorem 1.1} Let $F$ satisfy {\rm (F)} and let $\alpha , \beta \in C_0^2
({\mathbb T}, {\mathbb R})$. Then the asymptotics of
$\lambda ^N_n$ are as follows:

at the {\rm left edge}: for $0 \leq n \leq 2[F(N)]$
   \[ \lambda ^N_n = - 2 + \frac{1}{4N^2} \lambda ^-_n + O(F(N)^2 N^{-3})
   \]
at the {\rm right edge}: for $0 \leq n \leq 2[F(N)]$
   \[ \lambda ^N_{2N-1-n} = 2 - \frac{1}{4N^2} \lambda ^+_n + O(F(N)^2 N^{-3})
   \]
in the {\rm bulk}: for $ n = 2\ell , 2 \ell - 1$ with $[F(N)] < \ell < N - [F(N)]$,
    \[ \lambda ^N_n = - 2 \cos \frac{\ell \pi }{N}
       + O(N^{-2}F(N)^{-1}) .
   \]
These estimates hold uniformly in $0 \leq n \leq 2N-1$ and uniformly on bounded subsets of
functions $\alpha , \beta $ in $C_0^2 ({\mathbb T}, {\mathbb R})$.
\end{guess}

\begin{remark}
\label{Remark2.1}
In the case where $F(N)= N^\eta, 0 < \eta < 1/2,$ the asymptotics of Theorem \ref{Theorem 1.1} read
as follows:
\[ 
\lambda ^N_n = - 2 + \frac{1}{4N^2} \lambda ^-_n + O( N^{-3+2\eta}), 
\quad \forall \, 0 \leq n \leq 2[N^\eta]
\]
\[ 
\lambda ^N_{2N-1-n} = 2 - \frac{1}{4N^2} \lambda ^+_n + O( N^{-3+2\eta}),
\quad \forall \, 0 \leq n \leq 2[N^\eta]
\]
\[ 
\lambda ^N_{2\ell}, \lambda^N_{2\ell -1} = - 2 \cos \frac{\ell \pi }{N} + O(N^{-2- \eta}),
\quad \forall \, [N^\eta] < n < N - [N^\eta] .
\]
\end{remark}

To prove Theorem \ref{Theorem 1.1} we use singular perturbation methods, more specifically methods from semiclassical approximation. Indeed it has been proved in \cite{BGPU} that  Jacobi matrices such as $Q_N^{\alpha, \beta}$ can be viewed as matrix representations of certain 
semiclassical Toeplitz operators $T_N$ in the framework of the geometric quantization of the 2d torus. Note that the Jacobi matrices $\jmq$ -- and hence the associated  Toeplitz operators --  are perturbations of size $\frac1{N^2}$ of the equilibrium matrix $Q(0_N,1_N)$ whose spectrum is $\{-2\cos\frac{l\pi}N,\ l=0, \dots, N\}$. Since $\cos\frac{l\pi}N-\cos\frac{(l-1)\pi}N=O(N^{-2})$ for $l\sim 1$ or $l\sim N$, the size of the perturbation is of the same order as the spacing between the unperturbed eigenvalues so that regular perturbation methods fail. Using semiclassical methods, we compute the asymptotics of the eigenvalues at the two edges of the spectrum by constructing semiclassical (Lagrangian)  quasimodes for $T_N,$ for which the two Hill operators appear in the transport equation associated to the construction. As customary for the quantization of compact symplectic manifolds, the Toeplitz operators $T_N$ act on a Hilbert space of dimension $2N$ with $N$ playing the role of t
 he inverse of an effective Planck constant. In the bulk, i.e. for $1\ll l\ll N,$ one has $\vert\cos\frac{l\pi}N-\cos\frac{(l-1)\pi}N\vert\gg N^{-2}$ and thus regular perturbation methods apply.
Finally let us mention that our method allows to obtain the full asymptotic expansion in $\frac1{N^2}$ of the entire spectrum -- see the discussion at the end of Section \ref{eigenvalues}. 

\medskip

As an application of Theorem \ref{Theorem 1.1} we derive asymptotics for the characteristic polynomial $\chi_N(\mu)$ of $\jmq$ as \Ninf.
Note that $\chi_N(\mu)$ gives rise to the spectral curve $\{(\mu,z)\in\mathbb C^2\vert z^2=\chi_N(\mu)\}$ which plays an important role in the theory of periodic Toda lattices. These asymptotics will be of great use in the subsequent work \cite{bkp2}. By Floquet theory, $\chi_N(\mu)$ can be expressed in terms of the discriminant associated to the difference equation ($k\in\mathbb Z$)
\be\label{14bis}
  a_{k-1}^Ny(k - 1) + b_k^Ny(k) + a_k^N y(k + 1) = \mu y(k).
   \ee
Indeed denote by $y_1^N$ and $y_2^N$ the fundamental solutions of \eqref{14bis}
determined by

\[ y_1^N(0,\mu ) = 1,\, \, y_1^N(1,\mu ) = 0 \,\,\mbox { and } \,\, y_2^N(0,\mu ) = 0, \,\, y_2^N(1,\mu )
      = 1.
      \]

The discriminant of \eqref{14bis} is then defined as the trace of the Floquet matrix associated to \eqref{14bis} and given by
\[ \Delta _N(\mu ) = y_1^N(N,\mu ) + y_2^N(N + 1, \mu ).
   \]
   In view of the Wronskian identity, $\mu$ is an eigenvalue of \lab\  [$\jmq$] iff $\Delta_N(\mu)-2=0$ [$\Delta^2_N(\mu) - 4 = 0\}$].
Hence up to a multiplicative constant, $\Delta_N^2-4$ and $\chi_N$ coincide. 
From the recursive formula for $y_2^N(k, \mu )$  one then sees (\cite{HK1}) that $\Delta_N^2(\mu)-4=q_N^{-2}\chi_N(\mu)$ where $q_N=\prod\limits_1^N(1+\frac 1{4N^2}\alpha(\frac nN))$.

Analogously, denote by $\Delta_\pm(\lambda) \equiv \Delta(\lambda,q_\pm)$ the discriminant of 
\be\label{14ter}
-y''(x,\lambda)+q_\pm(x)y(x,\lambda)=\lambda y(x,\lambda)
\ee
defined as the trace of the Floquet operator associated to \eqref{14ter},
\[
\Delta_\pm(\lambda)=y_1^\pm(1/2,\lambda)+(y_2^\pm)'(1/2,\lambda)
\]
where $y_1^\pm(x,\lambda)$ and $(y_2^\pm)'(x,\lambda)$ are the fundamental solutions of \eqref{14ter} defined by
\[
y_1^\pm(0,\lambda)=1,\ (y_1^\pm(0,\lambda))'=0\,\,\,\mbox{ and }\,\,\,
y_2^\pm(0,\lambda)=0,\ (y_2^\pm(0,\lambda))'=1.
\]
Similarly as in the case of the Toda lattice, $\lambda$ is a periodic eigenvalue of $H_\pm$ on the interval $[0,1]$ iff $\Delta_\pm^2((\lambda)-4=0$.
Note that $\Delta_\pm^2((\lambda)-4$ is an entire function and can be viewed as a regularized determinant of $H_\pm$, referred to as characteristic function of \eqref{14ter}. It leads to the spectral curves $\{(\lambda,z)\in\mathbb C^2\vert \, z^2=\Delta_\pm^2(\lambda) - 4\}$ which play an important role in the theory of the KdV equation.
We will state our result on the asymptotics of $\chi_N$ in terms of the discriminant $\Delta_N$.
\medskip
With $M = [F(N)]$ and $F$  as before, let $\Lambda ^{\pm ,M}  \equiv \Lambda^{\pm , M}_2$ be the box
   \[ \Lambda ^{\pm ,M} := [\lambda ^\pm _0 - 2, \lambda ^\pm_{2[F(M)]} + 2] + i[-2,2]
   \]
and choose $N_0 \in {\mathbb Z}_{\geq 1}$ so that
   \be\label{81bis} \lambda ^\pm_{2k+1} - \lambda ^\pm _{2k} \geq 6 \quad \forall \
   k \geq F(F(N_0))
   \ee
By the Counting Lemma for periodic eigenvalues (cf e.g. \cite{KP}), $N_0$ can be chosen 
uniformly for
bounded subsets of function $\alpha , \beta $ in $C_0^2({\mathbb T}, {\mathbb R})$.

\begin{guess}
\label{Theorem9.1} Let $F$ satisfy {\rm (F)}, $M=[F(N)]$ with $N \ge N_0$,
and  $\alpha , \beta \in C_0^2({\mathbb T}, {\mathbb R})$.  Then uniformly for $\lambda$ in $\Lambda^{-,M}$
   \begin{equation}
    \label{9.1} \Delta_N \big(-2+\frac{1}{4N^2}\lambda\big) = (-1)^N\Delta_-(\lambda) + O
                \left(\frac{F(M)^2}{M} \right).
    \end{equation}
 Similarly, uniformly for $\lambda$ in $\Lambda^{+,M}$
  \begin{equation}
  \label{9.2} \Delta_N (2-\frac{1}{4N^2}\lambda) = \Delta_+(\lambda) + O\left(\frac{F(M)^2}
              {M}\right).
  \end{equation}
 These estimates hold uniformly on bounded subsets of functions $\alpha$ and $\beta $ in 
$C_0^2({\mathbb T}, {\mathbb R})$.
\end{guess}
\begin{remark}
In the case where $F(N)= N^\eta, 0 < \eta < 1/2,$ the asymptotics of Theorem \ref{Theorem9.1} 
read as follows:
\[
\Delta_N \big(-2+\frac{1}{4N^2}\lambda\big) = (-1)^N\Delta_-(\lambda) + O\left( N^{-\eta(1-\eta)} \right),
\quad \forall \, \lambda \in \Lambda^{-,N^{\eta}} ,
\]
\[
\Delta_N (2-\frac{1}{4N^2}\lambda) = \Delta_+(\lambda) + O\left( N^{-\eta(1-\eta)} \right),
\quad \forall \, \lambda \in \Lambda^{+,N^{\eta}} .
\]
\end{remark}

We remark that we did not aim at getting the maximal range of the
$\lambda $'s for which \eqref{9.1} and \eqref{9.2} hold. Moreover,
although we didn't state such a result here, our method allows to
compute the asymptotics of the discriminant for $\lambda$ in the bulk
region as well.  In the companion paper \cite{bkp2}, the results of
Theorems \ref{Theorem 1.1} and \ref{Theorem9.1} are used as an
important ingredient for computing the asymptotics of frequencies and
actions of Toda lattices in terms of the corresponding quantities of
the KdV equation.

\medskip

{\em Organisation of paper:} The proof of Theorem~\ref{Theorem 1.1} relies on the construction of quasimodes for the Jacobi matrices $\jm$. This construction is done in the framework of the geometric quantization of the torus (Section \ref{quantgeom}). The matrices $\jm$ are shown to be the matrix representation of Toeplitz operators (Section \ref{toplitz}), whose action on a certain type of Lagrangian states is described in detail in Theorem \ref{Proposition 3.5}.
Proposition \ref{Lemma 4.1} in Section \ref{quasilemma} states an abstract result on the construction of quasimodes that we couldn't find in the literature  and which is  crucial for the proof of Theorem \ref{Theorem 1.1}. The two cases corresponding to the bulk and the edges of the spectrum are treated in Section \ref{quasibulk} and Section \ref{quasiedges} respectively. The proof of Theorem \ref{Theorem 1.1} is summarized in Section \ref{eigenvalues}.
In Section \ref{discriminant} we first compute the asymptotics of the Casimir functionals of the Toda
lattice (Proposition \ref{jacobbb}) and then, using Theorem~\ref{Theorem 1.1},  obtain the asymptotics 
of the  discriminant of $\jm$ in terms of the discriminants of $H_{\pm}$, stated in Theorem \ref{Theorem9.1}. In addition, we apply Theorem \ref{Theorem9.1} to prove similar asymptotics for the derivatives of $\Delta_N(\mu)$ and to derive aymptotics of the zeroes of 
$\partial_\mu\Delta_N(\mu)$ at the two edges.

\medskip

{\em Methods:} The methodology used in this paper, based on the geometric quantization of the torus, is strongly inspired by \cite{BGPU}. In that paper, the authors consider the large $N$ asymptotics of Toda lattices, both for Dirichlet and periodic boundary conditions, in the case where the $a_n$'s and $b_n$'s are given by the discretization of regular functions,  i.e. the coupling parameter $\epsilon$ equals $1$, and  they derive the limiting PDE.

\medskip

Finally we mention that this work has been announced in \cite{BTP}.

\medskip

{\em Acknowledgments:} The authors would like to thank  the University of Milan, 
the Swiss National Science Foundation, the University of Z\"urich, the CNRS and the \'Ecole polytechnique
 for financial support during the elaboration of this work.


\section{Geometric quantization of $\mathbb T^2$}
\label{quantgeom}

 The geometric quantization of the two dimensional torus (resp.  sphere)   and the underlying so-called Toeplitz operators theory has been  shown in \cite{BGPU} to be a natural set-up for studying the large $N$ limit of the Toda lattice with periodic (resp. Dirichlet) boundary conditions. Although most of the computations in the present papers are going to be carried out from scratch, we recall in this section the basic facts concerning Toeplitz operators.


Consider the standard $2$-dimensional torus ${\mathbb T}^2 = {\mathbb R}^2
/ {\mathbb Z}^2$, identified with ${\mathbb C} / ({\mathbb Z} + i {\mathbb Z})$,
with canonical symplectic form $\omega = dx \wedge dy$ and Planck constant $(4\pi N)^{-1}$.
Let $E \rightarrow {\mathbb T}^2$ be a holomorphic line bundle with connection
$\nabla = d - 2\pi i x dy$ and denote by ${\kappa }$ the curvature form,
 ${\kappa } = d(-2\pi i x dy) .$
Then $\frac{i}{2\pi } \kappa = \omega $. In particular, the
Chern class of $E$, given by the cohomology class $\left[ \frac{i}{2\pi } {\kappa }
\right]$,
satisfies $\left[ \frac{i}{2\pi } {\kappa } \right] = [\omega ]$. Denote by
$({\mathcal H}^\sim _{2N}, \langle \cdot , \cdot \rangle _\sim )$ the Hilbert
space whose elements are holomorphic sections ${\mathbb T}^2 \rightarrow E^{\otimes
2N}$, viewed as entire functions $f : {\mathbb C} \rightarrow {\mathbb C}$ satisfying
   \[ f(z + m + in) =  e^{2N\pi [z(m-in) + \frac{1}{2}(m^2 + n^2)]} f(z)
   \qquad \forall (m,n) \in {\mathbb Z}^2, \,\,z \in {\mathbb C}. 
\]
The inner product is given by
 $\langle f,g \rangle _\sim = \int _{[0,1]^2} f(z) \overline{g(z)} e^{-2N\pi |z|^2} dxdy .$
We identify ${\mathcal H}^\sim _{2N}$ (isometrically) with the space ${\mathcal H}
_{2N}$ of theta functions of order $2N$, i.e. entire functions $f : {\mathbb C}
\rightarrow {\mathbb C}$, satisfying
   \[ f(z + m + in) = e^{2N\pi (n^2 - 2inz)} f(z) \qquad \forall (m,n) \in {\mathbb Z}^2, \,\,z \in {\mathbb C}
   \]
with inner product
   $ \langle f,g \rangle = \int _{[0,1]^2} f(z) \overline{g(z)} e^{-4\pi Ny^2} dxdy .$
   For $0 \leq j \leq 2N - 1$, let
   \begin{equation}
    \label{1}   \theta ^N_j(z) = (4N)^{1/4} e^{-\pi j^2 / 2N} \sum _{n \in {\mathbb Z}}
                e^{-\pi (2N n^2 + 2jn)} e^{2\pi iz(j +2Nn)} .
   \end{equation}
One verifies that $(\theta ^N_j)_{0 \leq j \leq 2N-1}$ is an orthonormal basis of ${\mathcal H}
_{2N}$. Observe that in contrast to the 'standard' case of the quantization of a cotangent bundle, the Hilbert space ${\mathcal H} _{2N}$ is finite dimensional. From a point of view of physics, this fact is justified by the Heisenberg uncertainty principle.
The Toeplitz quantization of a function $F : {\mathbb T}^2 \rightarrow {\mathbb R}$
is given by the sequence of operators
$ T^N_F : {\mathcal H}_{2N} \rightarrow {\mathcal H}_{2N} , \ f \mapsto P_N(fF)$
where $P_N : (L^2([0,1]^2, e^{-4N\pi y^2} dxdy) \rightarrow {\mathcal H}_{2N}$ denotes the orthogonal projector,
   \begin{align*} 
                  &(P_Nf)(z) = \sum ^{2N-1}_{j=0} \langle f, \theta ^N_j \rangle \theta
                     ^N_j(z) .
   \end{align*}
More generally, a Toeplitz operator is a sequence of operators $(T^N)_{N \geq 1}$ where
for $N \geq 1, T^N : {\mathcal H}_{2N} \rightarrow {\mathcal H}_{2N}$ is an operator of
the form
   $ T^N \sim \sum ^\infty _{j = 0} N^{-j} T^N_{S_j}.
   $
The function $S_0 : {\mathbb T}^2 \rightarrow {\mathbb R}$ is referred to as principal symbol. 


\section{Jacobi matrices as Toeplitz operators}
\label{toplitz}

In this section we study Jacobi matrices with entries defined in terms of discretizations of functions of a certain regularity, from a 'Toeplitz operator' perspective. In particular we show how their action on 
certain families of elements in the Hilbert spaces $\mathcal H_{2N}$, referred to as Lagrangian states, 
 can be explicitly described in terms of differential operators -- see Theorem \ref{Proposition 3.5} below. This theorem is key for our construction of quasimodes.

For $\alpha , \beta $ in $C^2_0({\mathbb T}, {\mathbb R})$ and $N \geq 3$,
denote by $T^{\alpha , \beta }_N$ the linear operator on ${\mathcal H}_{2N}$ whose
representation with respect to the basis $[\theta ^N_{2N-1}, \ldots , \theta ^N_0]$ is
$Q^{\alpha , \beta }_N$. As an example, consider the operator $T^{0,0}_N$. To
study its properties let us begin by recording the following elementary result.
\begin{lemma}
\label{Lemma 3.1}
$\big( (\mp1)^n\big)_{0 \leq n \leq 2N-1}$ is an eigenvector of $Q^{0,0}_N$
corresponding to the eigenvalue $\mp2$, and, for any $1 \leq \ell \leq N - 1$, the vectors
  $ \big( e^{i\pi (N-\ell )n/N} \big) _{0 \leq n \leq 2N-1} $ and
      $\big( e^{i\pi (N+\ell )n/N} \big) _{0 \leq n \leq 2N-1}
  $
are eigenvectors of $Q^{0,0}_N$ corresponding to the eigenvalue $-2 \cos \frac{\ell
\pi }{N}$.
They form an orthogonal basis in ${\mathbb C}^{2N}$.
\end{lemma}

From Lemma~\ref{Lemma 3.1} it follows that $\psi ^{N,k}(z), \ 0 \leq k \leq 2N-1,$
is an orthonormal basis of ${\mathcal H}_{2N}$ of eigenfunctions of $T^{0,0}_N$,
$\,T^{0,0}_N \psi ^{N,k} = 2 \cos \frac{k\pi }{N}  \psi ^{N,k}$,  where
   \begin{align}
   \begin{split} \psi ^{N,k}(z) &= (2N)^{-1/2} \sum ^{2N}_{n=1} e^{i \pi \frac{k n}
                     {N}} \theta ^N_{2N-n}(z)  = (2N)^{-1/2} \sum ^{2N-1}_{n=0} e^{-i \pi \frac{k n}
                     {N}} \theta ^N_{n}(z)\label{psiN} .
   \end{split}
   \end{align}
 Alternatively, $\psi ^{N,k}$ can be expressed with the help of the kernel
   \begin{equation}
   \label{2} \rho _N(z,w) = \sum ^{2N-1}_{j=0} \theta ^N_j(z) \overline {\theta ^N_j(w)} .
   \end{equation}

\begin{lemma}
\label{Lemma 3.2} For any $0 \leq k \leq 2 N - 1$,
   \[ \psi ^{N,k}(z) = (4N)^{-1/4} \int ^1_0 \rho _N \big( z, k/2N + is\big)
      e^{-2N\pi s^2} ds.
   \]
\end{lemma}

{\it Proof: } In view of \eqref{1} and \eqref{2}, the claimed identity follows easily from the identity
   $ \sum _{n \in {\mathbb Z}} \int ^1_0 e^{-2\pi N(n + s + \frac{j}{2N})^2} ds = \int ^\infty _{-\infty }
      e^{-(\sqrt{2\pi N}x)^2} dx = (2N)^{-1/2}.
   $\hspace*{\fill }$\square$

\medskip

It is useful to introduce for an arbitrary real or complex valued function $f \in L^2
({\mathbb T})$ and $0 \leq k \leq 2N - 1$,
   \begin{equation}
   \label{4} \psi ^{N,k}_f(z) = (4N)^{-1/4} \int ^1_0 f(s) \rho _N \big( z,k/2N
             + is \big) e^{-2N\pi s^2} ds .
   \end{equation}
It is an element in ${\mathcal H}_{2N}$. For $f \in L^2({\mathbb T})$, denote by $\hat
f_n$ the $n$'th Fourier coefficient of $f$,
   $\hat f_n = \int ^1_0 f(x) e^{-i2\pi nx} dx
   $
and by $\| f\| _\ell $ the $\ell $-th Sobolev norm
   $ \| f\| _\ell = \Big( \sum _{n \in {\mathbb Z}} |\hat f_n|^2 (1 + |n|)^{2\ell }
      \Big) ^{1/2} .
   $
Further denote by $\| f\| _{C^\ell }$ the following norm of $f \in C^\ell ({\mathbb T})$:
   $\| f\| _{C^\ell } = \sup\limits _{0 \leq x \leq 1} \sum ^\ell _{j=0} |\partial ^j_x
      f(x)| .$

\smallskip

\begin{lemma}
\label{Lemma 3.3} For $f, g \in L^2({\mathbb T})$ and $0 \leq k, \ell \leq 2N - 1$

\begin{description}
\setlength{\labelwidth}{9mm}
\setlength{\labelsep}{1.4mm}
\setlength{\itemindent}{0mm}

\item[{\rm (i)}] $\psi ^{N,k}_f(z) = \frac{1}{\sqrt{2N}} \sum ^{2N-1}_{j=0} \theta ^N
_j(z) e^{-i\pi kj/N} \sum _{m \in {\mathbb Z}} \hat f_m e^{-\pi m^2/2N} e^{-i\pi mj/N}$.

Alternatively, with $\Delta = -d^2/dx^2$,
 \be\label{heatt} 
   \psi ^{N,k}_f (z) = \frac{1}{\sqrt{2N}} \sum ^{2N-1}_{j=0}  \big( e^{-\Delta /8\pi N} f \big) 
    \big( -j/2N    \big)e^{-i \pi kj/N} \theta ^N_j(z). 
 \ee 
\item[{\rm (ii)}] $\langle \psi ^{N,k}_f, \psi ^{N,\ell }_g \rangle = \sum _{n \in {\mathbb Z}}
\hat f_n \overline {\hat g_{n + k - \ell }} e^{-\pi n^2/2N} e^{-\pi (n + k - \ell )^2/2N}$.

\item[{\rm (iii)}] 
   $ \big\arrowvert\langle \psi ^{N,k}_f, \psi ^{N,k}_g \rangle - \langle f,g \rangle
      \big\arrowvert \leq \frac{1}{4\pi N} \| f'\| _0 \cdot \|g'\|_0 $\quad $\forall \,f,g \in H^1({\mathbb T})$.

\item[{\rm (iv)}] The linear maps $L^2({\mathbb T}) \rightarrow {\mathcal
H}_{2N}, f \mapsto \psi ^{N,k}_f,$ are bounded, 
   $\|\psi ^{N,k}_f \| \leq \| f\|_0.$
\end{description}
\end{lemma}

{\it Proof: } (i) By the definitions of $\psi ^{N,k}_f$ and $\rho _N$
   \[ \psi ^{N,k}_f(z) = (4N)^{-1/4} \sum ^{2N-1}_{j=0} \theta ^N_j(z) \int ^1_0 f(s)
      \ \overline{ \theta ^N_j ( k/2N + is) } \  e^{-2N \pi s^2}
      ds .
   \]
Using the definition \eqref{1} of $\theta _j^N$ one gets
   \begin{align} \label{6} &(4N)^{-1/4} \int ^1_0 f(s) \overline{\theta ^N_j (\frac k{2N} +
                     is)} e^{-2N\pi s^2} dx
                     = e^{-i\pi\frac{kj}N} \sum _{n \in {\mathbb Z}} \int ^1_0 e^{-2\pi N (n + s +\frac{ j}{2N})^2}
                     f(s) ds.
   \end{align}

As $ \sum _{n \in {\mathbb Z}} e^{-2\pi N(n + s + j/2N)^2}= \frac{1}{\sqrt{2N}} \sum _{n\in {\mathbb Z}} e^{2\pi n(s + j/2N)} e^{-\pi n^2/2N}$ (Poisson summation formula) it then follows that

   \begin{align*}
&{(4N)^{-1/4}} \int ^1_0 f(s) \overline{\theta ^N_j \big( k/2N +
                     is \big)} e^{-2N\pi s^2} ds = \frac{1}{\sqrt{2N}} \sum _{n \in {\mathbb Z}} e^{i\pi \frac{(n-k)j}N}
                     e^{-\pi \frac{n^2}{2N}} \hat f_{-n},
   \end{align*}
yielding the claimed formula. To verify the alternative formula, note that by \eqref{6} and the fact that $f$ is periodic, one has
   
   \begin{align*} \psi ^{N,k}_f(z) &= \sum ^{2N-1}_{j=0} \theta ^N_j(z) e^{-i\pi kj/N}
                     \sum _{n \in {\mathbb Z}} \int ^1_0 e^{-2\pi N(n + s + j/2N)^2}
                     f(s)ds \\
                  &= \sum ^{2N-1}_{j=0} \theta ^N_j(z) e^{-i\pi kj/N} \int ^\infty _{-
                     \infty } e^{-2\pi N(y + j/2N)^2} f(y)dy.
   \end{align*}
 Evaluating the heat flow for the initial data $f$ at $x= -\frac{ j}{2N}$ and time $t= \frac{1}{4}\frac{1}{2\pi N}$  we obtain the claimed identity \eqref{heatt}.

\smallskip

(ii) By the definition of $\rho _N$ and the fact that $(\theta ^N_j)_{0 \leq j \leq 2N-1}$ is an orthonormal basis of $\mathcal H _{2N}$, we have that
   \begin{align*} &\int ^1_0 \rho _N \big( x + iy , k/2N + is \big) \rho _N
                     \big( x + iy, \ell /2N + it \big) e^{-4\pi Ny ^2} dx dy \\
                  & =\sum _{j,n} \overline{\theta ^N_j \big( k/2N + is \big) }
                     \theta ^N_n \big( \ell / 2N + it \big) \langle \theta ^N
                     _j, \theta ^N_n \rangle = \sum ^{2N-1}_{j=0} \overline{\theta ^N_j \big( \frac k{2N} + is \big) }
                     \theta ^N_n \big( \frac\ell {2N} + it \big).
   \end{align*}
Hence $\langle \psi ^{N,k}_f, \psi ^{N,\ell }_g \rangle $ is equal to
   \begin{align*} &\sum ^{2N-1}_{j=0} \frac1{\sqrt{4N}} \int \limits^1_0\overline{\theta ^N_j \big(
                     k/2N + is \big) } f(s) e^{-2N\pi s^2} ds  \int \limits^1_0 \theta ^N_j \big( k/2N + it \big)
                     \overline {g(t)} e^{-2N\pi t^2}dt.
   \end{align*}
By (i) and the fact that $\sum ^{2N-1}_{j=0} e^{i\pi (\ell - n - k + m)j/N} =  2N\delta_{\ell-n,k-m}$ we  get that
   \begin{align*} \langle \psi ^{N,k}_f, \psi ^{N,\ell }_g \rangle 
                 &= \sum _{m \in {\mathbb Z}} \hat f_{-m} \overline {\hat g_{-m + k - \ell }}
                    e^{-\pi m^2/2N} e^{-\pi (-m + k - \ell )^2/2N}.
   \end{align*}

\smallskip

(iii) From item (ii) it follows that for $k = \ell $,
   \begin{align*} \langle \psi ^{N,k}_f, \psi ^{N,k}_g \rangle &= \sum _{n \in {\mathbb
                     Z}} \hat f_n \overline {\hat g_n} e^{-\pi n^2/N} = \langle f, g\rangle - \sum _{n \in {\mathbb Z}} \hat f_n \overline{\hat g
                     _n} (1 - e^{-\pi n^2/N}) .
   \end{align*}
As $0 \leq 1 - e^{-\pi n^2/N} \leq \pi n^2 / N$ one has by the Cauchy-Schwarz inequality
   \begin{align*} \big\arrowvert \sum _{n \in {\mathbb Z}} \hat f_n \hat g_n (1 - e^{-\pi n^2 /
                     N}) \big\arrowvert &\leq \frac{1}{4\pi N} \sum _{n \in {\mathbb Z}}
                     |\hat f_n| |\hat g_n| (2\pi n)^2 \leq \frac{1}{4\pi N} \| f'\| _0 \| g'\| _0
   \end{align*}
and the claimed estimate follows.
\smallskip

(iv) By (ii) (for $f = g$ and $k = \ell $)
   $\| \psi ^{N,k}_f \| ^2 = \sum _{n \in {\mathbb Z}} |\hat f_n| ^2 e^{-\pi n ^2 / N}
      \leq \| f\| ^2_0 .
      $
\hspace*{\fill }$\square $

\medskip

Next we describe how $T^{\alpha , \beta  }_N$ acts on $\psi ^{N,k}_f$. This result will be
an important ingredient to obtain the asymptotics of the eigenvalues of $Q^{\alpha ,
\beta }_N$ at the edges. For $f \in L^2({\mathbb T})$ and $0 \leq \ell \leq 2N - 1$,
introduce, with $\alpha _2(x):= \alpha (2x)$,
   \begin{align}\label{nouveau2} D^{\alpha , \beta }_\ell (f):= &2 \cos \big( \frac{\ell \pi} N -
                     \frac i{2N} \partial _x \big) f + \frac1{4N^2} \big( \beta _2(x) + 2\alpha _2(x) \cos \big(
                     \frac{\ell \pi }N - \frac i{2N} \partial _x \big) \big) f .
   \end{align}
This expression is to be understood in the sense of functional calculus. More precisely,
$\cos \big( \frac{\ell \pi }{N} - \frac{i}{2N} \partial _x \big) $ is viewed as a
multiplier operator in Fourier space
   \[ \cos \big( \ell \pi /N - i/2N \partial _x \big) f = \sum _{n \in
      {\mathbb Z}} \hat f_n \cos \big( \ell \pi /N + 2\pi n/2N \big)
      e^{i2\pi nx} .
   \]

\begin{guess}
\label{Proposition 3.5} For any $f \in C^2$ and $0 \leq \ell \leq 2N - 1$
   \[ \| T^{\alpha , \beta }_N \psi ^{N,\ell }_f - \psi ^{N,\ell }_{D^{\alpha , \beta }_\ell
      (f)}\| \leq \frac{1}{N^3} K _{\alpha , \beta } \| f\|_{C^2}
   \mbox{ with }K _{\alpha , \beta }:= \| \alpha \| _{C^2} + \| \beta \| _{C^2} + 1.\]
\end{guess}

To prove Theorem~\ref{Proposition 3.5} we first need to establish some auxiliary results.
First note that  $Q^{\alpha , \beta }_N = Q^{0, \beta }_N + \big( Q^{0,\alpha }_N -
                     Q^{0,0}_N \big) Q^+_N + Q^-_N \big( Q^{0,\alpha }_N - Q^{0,0}_N \big)$
where
   \[ Q^+_N = \begin{pmatrix} 0 &1 & &0 \\ 0 &0 & &\vdots \\ \vdots & \vdots
      &\ldots &1 \\ 1 & 0 & & 0  \end{pmatrix}
   \]
and $Q^-_N$ is the transpose of $Q^+_N$. Denote by $T^\pm _N$ the operator on
${\mathcal H}_{2N}$ whose matrix representation with respect to the basis $[\theta ^N_{2N-1},
\ldots , \theta ^N_0]$ are $Q ^\pm_N$. Notice that $T^\pm _N$ are isometries
as $Q^\pm _N$ are the matrix representations of permutations. Further $T^{0,0}
_N = T^+_N + T^-_N$ as $Q^{0,0}_N = Q^+_N + Q^-_N$. For any $f \in L^2({\mathbb T})$
and $0 \leq \ell \leq 2N - 1$, define
   \[ D^\pm _\ell (f) = \exp \big( \pm i \frac{2\pi \ell - i \partial/\partial x}{2N} \big)f
      \,\, \mbox{ and } \,\,D^{0,0}_\ell = D^+_\ell + D^-_\ell .
   \]

\begin{lemma}
\label{Lemma 3.6} For any $f \in L^2({\mathbb T})$ and $0 \leq \ell \leq 2N-1$,

\[
T^\pm_N \psi ^{N,\ell }_f = \psi ^{N,\ell }_{D^\pm _\ell (f)}\quad\mbox{ and } \quad T^{0,0}_N \psi ^{N,\ell }_f = \psi ^{N,\ell }_{D^{0,0} _\ell (f)}.
\]
\end{lemma}

\smallskip

{\it Proof:}  Since $\psi ^{N,\ell }_f$ is linear in $f$ it suffices to verify the claimed identities for
$f(x) =e_k(x):= e^{i2\pi kx}$.
By Lemma~\ref{Lemma 3.3} (i) and the fact that $(e^{i\pi nj/N})_{0 \leq j \leq 2N-1}$ is
an eigenvector of $Q^\pm_N$ with eigenvalue $e^{\pm i\pi n/N}$ one has
\[ T^\pm _N \psi^{N,\ell }_{e _k} = e^{\pm i \pi (k + \ell )/N} \psi ^{N,\ell }
      _{e_k} = \psi ^{N,\ell }_{e^{\pm i \pi (k + \ell )/N}e_k}
   \]
As $e^{\pm i \pi (k + \ell )/N} e^{i2\pi kx} = e^{\pm i\pi \ell/N} e
                     ^{\pm \frac {-i\partial /\partial x}{2N}} \big( e^{i2\pi k x}\big) 
                  = D^\pm _\ell \big( e^{i2\pi kx} \big)$
 the claimed identities  follow.
 \hspace*{\fill }$\square$

\bigskip

The key result used in the proof of Theorem~\ref{Proposition 3.5} is the following one.
\begin{lemma}
\label{Lemma 3.7} Let $f,g \in C^2({\mathbb T})$. Then with $f_2(x):= f(2x)$, one has
   \[ \\| (T^{0,f}_N - T^{0,0}_N) \psi ^{N,k}_g - \frac{1}{4N^2} \psi^{N,k}_{gf_2}\|
      \leq \frac{1}{32\pi N^3} \big( \| (gf_2)'' \| _{C^0} + \| f \| _{C^0}
      \| g''\| _{C^0} \big).
   \]
\end{lemma}

{\it Proof:} Note that $Q^{0,f}_N - Q^{0,0}_N$ is a diagonal
matrix with entries $(2N)^{-2} f(j/N)$, $1 \leq j \leq 2N$.
By the definition of $T^{0,f}_N$, it then follows that
   \[ \big( T^{0,f}_N - T^{0,0}_N \big) \theta ^N_j = (2N)^{-2} f \big( (2N-j)/
      N \big) \theta ^N_j = (2N)^{-2} f \big( - j/N \big) \theta ^N_j .
   \]
Hence by Lemma \ref{Lemma 3.3} (i)
\[4N^2 \big( T^{0,f}_N - T^{0,0}_N \big) \psi ^{N,k}_g(z) = \frac{1}{\sqrt{2N}}
      \sum ^{2N-1}_{j=0} f_2 ( - \frac j{2N} ) \theta ^N_j(z)
                  e^{-i\pi kj/N}\big(e^{-\frac\Delta{8\pi N}}g \big) ( - \frac j{2N}) .
\]
Furthermore, one has
   $\psi ^{N,k}_{gf_2} (z) = \frac{1}{\sqrt{2N}} \sum ^{2N-1}_{j=0} \theta ^N_j(z) e^{-i
      \pi kj/N} e^{-\Delta /8\pi N} (gf_2) \big( -j/2N \big) .
   $
As $(\theta _j)_{0 \leq j \leq 2N-1}$ is an orthonormal basis of ${\mathcal H}_{2N}$, it
then follows
   \begin{align*} &\| \psi ^{N,k}_{f_2g} - 4N^2 (T^{0,f}_N - T^{0,0}_N) \psi ^{N,k}_g
                     \| ^2 \leq \frac{1}{2N} \sum ^{2N-1}_{j=1} \big\arrowvert [ e^{-\Delta / 8\pi
                     N}, M_{f_2}] g \big( - \frac j{2N}\big) \big\arrowvert ^2
   \end{align*}
where $M_{f_2}$ denotes the operator on $L^2({\mathbb T})$ of multiplication by $f_2$ and
$[\cdot , \cdot ]$ is the commutator of operators. Hence
   \[ \| \psi ^{N,k}_{f_2 g} - 4N^2 (T^{0,f}_N - T^{0,0}_N) \psi ^{N,k}_g \| \leq \sup _{0
      \leq x \leq 1} \big\arrowvert [e^{-\Delta /8\pi N}, M_{f_2}] g(x) \big\arrowvert .
   \]
We estimate the latter expression using
   $e^{-\Delta t} = Id - \Delta \int ^t_0 e^{-\Delta s} ds$,
    \[ f_2(x)(e^{-\Delta / 8\pi N} g) (x) = f_2(x)g(x) - f_2(x) \left( \int ^{(8\pi N)^{-1}}_0
      e^{-\Delta s} ds \, \, \Delta g \right) (x) .
   \]
Using the formula of the heat kernel on $\mathbb R$ and the identity $\int ^\infty _{-\infty } e^{-(x-y)^2/4s} dy = \sqrt{4\pi s}$ one gets
   \begin{align*} \big\arrowvert \big( \int ^{(8\pi N)^{-1}}_0 e^{-\Delta s} ds\,\, \Delta g\big)
                     (x) \big\arrowvert \leq \| g'' \| _{C^0} (8\pi N)^{-1} \quad \mbox{ and }
   \end{align*}
   \[ \big\arrowvert \int ^{(8\pi N)^{-1}}_0 e^{-\Delta s} ds \Delta (gf_2) (x)
      \big\arrowvert \leq \| (gf_2 )''\| _{C^0} \ (8\pi N)^{-1} .
   \]
Combining these estimates yields the claimed estimate. \hspace*{\fill}{$\square $}

\bigskip

Finally, for the proof of Theorem~\ref{Proposition 3.5} we will also need the following lemma.
\begin{lemma}
\label{Lemma 3.8} For $f \in C^1({\mathbb T})$, denote by $M_{f_2}$ the multiplication
operator on $L^2({\mathbb T})$ by $f_2(x):= f(2x)$. Then the  operator $[D^\pm
_k, M_{f_2}]$ on $L^2({\mathbb T})$ satisfies
   \[ [D^\pm _k, M_{f_2}] = \big( f(2x) - f(2x \pm \frac{1}{N}) \big) D^\pm _k
   \quad \mbox{and} \quad
    \| [D^\pm _k, M_{f_2}] \| _{L^2 \rightarrow L^2} \leq \frac{1}{N} \| f'\| _{C^0}.\]
   Moreover
   \[\big\| [T^{0,f}_N- T^{0,0}_N, T^\pm _N ]\big\| _{{\mathcal H}_{2N} \rightarrow {\mathcal H}
      _{2N}}\leq \frac{1}{N} \| f'\| _{C^0}.
   \]
\end{lemma}
{\it Proof: } Recalling  that  $D^\pm _k = e^{\pm i \frac{2\pi k - i
\partial / \partial x}{2N}} $, it is straightforward to verify that the values of the two operators coincide for any $g = e^{i2\pi nx}$. The claimed identity then follows by linearity. The claimed bound of the operator norm of $[D^\pm _k, M_{f_2}]$ then follows from the unitarity of $D^\pm _k$. 
The second estimate is proved using the matrix representation of the operators involved. 
Recall that $Q^{0,f}_N - Q^{0,0}
_N = {\rm diag} \big( (2N)^{-2}f(j/N)_{1 \leq j \leq
      2N} \big)$.
Thus $[(Q^{0,f}_N - Q^{0,0}_N), Q^+_N]$ is the $2N \times 2N$ matrix
 \[ \begin{pmatrix} 0 & \gamma _1 & &  0 \\ 0&0&\ldots &\vdots \\ \vdots &\vdots &
      & \gamma _{2N-1} \\ \gamma _{2N}&0 &&0 \end{pmatrix}
 \quad 
\mbox{with } \quad 
\gamma _i = f\big( \frac{i}{N} \big) - f \big( \frac{i + 1} {N} \big). 
\]
Hence
   \[ \big\| [ (Q^{0,f}_N - Q^{0,0}_N), Q^+_N] \big\| _{{\mathbb R}^{2N}
                     \rightarrow {\mathbb R}^{2N}} = \sup _{1 \leq i \leq 2N} \big\arrowvert
                     f\big( \frac{i}{N} \big) - f \big( \frac{i + 1}{N} \big)
                     \big\arrowvert \leq \frac{1}{N} \| f' \| _{C^0} .
   \]
As $Q^-_N(Q^{0,f}_N - Q^{0,0}_N)$ is the transpose of $(Q^{0,f}_N - Q^{0,0}_N)Q^+_N$, the
same estimate holds for $[(Q^{0,f}_N - Q^{0,0}_N), Q^-_N]$.
\hspace*{\fill}$\square $

\medskip

{\it Proof of Theorem~\ref{Proposition 3.5}:} We write $T^{\alpha , \beta }_N$ as a sum
of operators
   \begin{align*} T^{\alpha , \beta }_N &= T^{0,0}_N + \big( T^{0,\beta }_N - T^{0,0}_N
                     \big) + \big( T^{0, \alpha }_N - T^{0,0}_N \big) T^+_N
                     + T^-_N \big( T^{0, \alpha }_N - T^{0,0}_N \big) \\
                  &= T^{0,0}_N + \big( T^{0,\beta }_N - T^{0,0}_N
                     \big) + T^{0,0}_N \big( T^{0, \alpha }_N - T^{0,0}_N \big)
                     + \left[ \big( T^{0, \alpha }_N - T^{0,0}_N \big) , T^+_N
                     \right] .
   \end{align*}
By Lemma~\ref{Lemma 3.6} and  Lemma~\ref{Lemma 3.7} we get, respectively,
$T^{0,0}_N \psi ^{N,\ell }_f = \psi ^{N,\ell }_{D^{0,0}_\ell (f)}$ and 
\[ \big\| \big( T^{0,\beta}_N - T^{0,0}_N \big) \psi ^{N,\ell }_f - \frac{1}
              {4N^2} \psi ^{N,\ell }_{\beta _2 f} \big\| \leq \frac{ 1}{8N^3} \| \beta \|
              _{C^2} \| f\| _{C^2} \quad \mbox { and }
\]
   \begin{align*} &T^{0,0}_N \big( T^{0,\alpha }_N - T^{0,0}_N \big) \psi ^{N,\ell }_f -
                     \frac{1}{4N^2} \psi ^{N,\ell }_{D^{0,0}_\ell (\alpha _2 f)} = T^{0,0}_N \big( (T^{0,\alpha }_N - T^{0,0}_N) \psi ^{N,\ell }_f -
                     \frac{1}{4N^2} \psi ^{N,\ell }_{\alpha _2 f} \big) .
   \end{align*}
As $T^{0,0}_N = T^+_N + T^-_N$ and $T^\pm _N$ are isometries it  follows from
Lemma~\ref{Lemma 3.7} that
   \[  \big\| T^{0,0}_N \big( T^{0,\alpha }_N - T^{0,0}_N \big) \psi ^{N,\ell }_f -
                     \frac{1}{4N^2} \psi ^{N,\ell }_{D^{0,0}_\ell (\alpha _2 f)} \big\|
                \leq \frac{1}{4N^3} \| \alpha \| _{C^2} \| f\| _{C^2} .
   \]
By Lemma~\ref{Lemma 3.8} and Lemma \ref{Lemma 3.3} (iv) it follows that
   \[  \big\| \big[ \big( T^{0,\alpha }_N - T^{0,0}_N\big), T^{+}_N \big]
                  \psi ^{N,\ell }_f \big\| \leq
                  \frac{\| \alpha ' \| _{C^0}}{4N^3} \| \psi ^{N,\ell }_f \|
                \leq \frac{\| \alpha ' \| _{C^0}}{4N^3}  \| f\| _{0} .
   \]
Finally, we need to estimate $\psi ^{N,\ell }_{\alpha _2 D^{0,0}_\ell f} - \psi ^{N,\ell }
_{D^{0,0}_\ell (\alpha _2 f)}$. As by
Lemma \ref{Lemma 3.3} (iv), the linear map $g \mapsto \psi ^{N,\ell }_g$ is bounded by $1$ 
on $L^2({\mathbb T})$,
it follows from Lemma~\ref{Lemma 3.8} that
   \begin{align*}
     \big\| \psi ^{N,\ell }_{\alpha _2 D^{0,0}_\ell f} - \psi ^{N,\ell }
                  _{D^{0,0}_\ell (\alpha _2 f)} \big\| &\leq \| \alpha _2 D^+_\ell f - D^+_\ell
                  (\alpha _2 f) \| _0 + \| \alpha _2 D^-_\ell f - D^-_\ell (\alpha _2 f)
                  \| _0 \\
               & \leq \frac{2}{N} \| \alpha ' \| _{C^0} \| f\| _{0} .
   \end{align*}
Taking into account the simple identity, $D^{\alpha , \beta }_\ell (f) = 
D^{0 , 0}_\ell (f) + \frac{1}{4N^2} \beta _2 f + \frac{1}{4N^2} \alpha _2 D^{0,0}_\ell (f)$,
 the obtained estimates imply the claimed one.\hspace*{\fill}$\square $


\section{Spectral results by quasimodes}
\label{quasilemma}

In this section we prove results on quasimodes used in the proof of Theorem~\ref{Theorem 1.1}.
Assume that ${\mathcal H}$ is a finite dimensional Hilbert space with inner product $\langle
\psi , \phi \rangle $ and induced norm $\| \psi \| = \langle \psi , \psi \rangle ^{1/2}$.
Further assume that $A : {\mathcal H} \rightarrow {\mathcal H}$ is a selfadjoint linear
operator.

\begin{proposition}
\label{Lemma 4.1} {\rm (i)} Assume that there exist $\psi \in {\mathcal H}$ with $\| \psi \|
= 1, \mu \in {\mathbb R}$ and $C > 0$ so that
   \begin{equation}
   \label{4.1} \| (A - \mu ) \psi \| \leq C .
   \end{equation}
Then there exists an eigenvalue $\lambda $ of $A$ so that $|\lambda - \mu | \leq C$.

{\rm (ii)} Assume that there exist two elements $\psi _\pm \in {\mathcal H}, \| \psi _\pm
\| = 1, \mu \in {\mathbb R}, 0 \leq \theta < 1$, and $C > 0$ so that
   \[ \| (A - \mu ) \psi _\pm \| \leq C \quad \mbox { and } \quad \big\arrowvert \langle \psi _+ ,
      \psi _- \rangle \big\arrowvert \leq \theta .
   \]
Then for any $D > 8C(1 - \theta )^{-1}$, there exist two eigenvalues $\lambda _\pm $ of
$A$ so that $|\lambda _\pm - \mu | \leq D$. If $\lambda _+ = \lambda _-$, then the
multiplicity of $\lambda _+$ is at least two.
\end{proposition}

\smallskip

{\it Proof:} (i) Denote by $(\lambda _j)_{j \in I}$ the eigenvalues of $A$ listed with
their multiplicities. As $A$ is selfadjoint ${\mathcal H}$ has an orthonormal basis of
eigenvectors, $(\psi _j)_{j \in I}$, where $\psi _j \in {\mathcal H}$ is an eigenvector
corresponding to the eigenvalue $\lambda _j$. Assume that for any $j \in I$,
   $| \lambda _j - \mu | > C .
   $
Then the vector $\psi = \sum _{j \in I} \langle \psi , \psi _j \rangle \psi _j$ satisfies
   \[ C^2 = C^2 \| \psi \| ^2 < \sum _{j \in I} \big\arrowvert \langle \psi , \psi _j \rangle
      \big\arrowvert ^2 (\lambda _j - \mu )^2 = \| (A - \mu ) \psi \| ^2 \leq C^2 ,
   \]
a contradiction. Hence the assumption is not true and (i) follows.

\smallskip

(ii) By item (i), there exists an eigenvalue $\lambda _{i_1}$ with $|\mu - \lambda _{i_1}|
\leq C$. Let us assume that $\lambda _{i_1}$ has multiplicity one and
   \begin{equation}
   \label{4.2} |\mu - \lambda | > D \quad \forall \lambda \in {\rm spec}(A) \backslash \{
               \lambda _{i_1}\} .
   \end{equation}
Then
   $P:= \frac{1}{2\pi } \int _K (z - A)^{-1} dz
   $
is the orthogonal projector of ${\mathcal H}$ onto the {\it one} dimensional eigenspace
of the eigenvalue $\lambda _{i_1}$ where $K$ denotes the counterclockwise oriented
circle of radius $D/2$ centered at $\lambda _{i_1}$. To estimate $P \psi _\pm $ note that
   $ \psi _\pm = (z - A)^{-1} (z - A) \psi _\pm = (z - A)^{-1} (z - \lambda _{i_1})
      \psi _\pm + (z - A)^{-1} r_\pm
   $
where $r_\pm = (\lambda _{i_1} - A)\psi _\pm $. Note that
   $\| r_\pm \| \leq \| (\mu - A) \psi _\pm \| + |\mu - \lambda _{i_1}|
        \leq 2C$
and
   \begin{equation}
   \label{A.2bis} (z - A)^{-1} \psi _\pm = (z - \lambda _{i_1})^{-1} \psi _\pm - (z - \lambda
                  _{i_1})^{-1}(z - A)^{-1} r_\pm .
   \end{equation}
Write $r_\pm$ as $r_\pm = Pr_\pm + (Id - P)r_\pm $ and use 
   $(z - A)^{-1} P r_\pm= (z - \lambda _{i_1})^{-1} P r_\pm
   $
to see that
   $ \frac{1}{2\pi i} \int _K(z - \lambda _{i_1})^{-1} (z - A)^{-1} P r_\pm dz = 0.$ Hence
   \begin{align*} &\frac{1}{2\pi i} \int _K(z - \lambda _{i_1})^{-1} (z- A)^{-1} r_\pm dz = \frac{1}{2
      \pi i} \int _K (z - \lambda _{i_1})^{-1}(z - A)^{-1} (Id - P) r_\pm dz.
   \end{align*}
By Cauchy's theorem we then get
   \begin{equation}
   \label{A.2ter} \frac{1}{2\pi i} \int _K (z - \lambda _{i_1})^{-1}(z - A)^{-1} r_\pm dz =
                  (\lambda _{i_1} - A)^{-1} (Id-P)r_\pm .
   \end{equation}
Hence integrating \eqref{A.2bis} along the contour $K$ one concludes from \eqref{A.2ter}
that
   \[ P \psi _\pm = \psi _\pm + (\lambda _{i_1} - A)^{-1} (Id - P) r_\pm .
   \]
By \eqref{4.1} - \eqref{4.2} we then have
   $ \| (\lambda _{i_1} - A)^{-1} (Id - P) r_\pm \| \leq D^{-1} 2C
   $
and thus, with $\eta := 2CD^{-1} < 1$, it follows that
   $ 0 \leq 1 - \| P \psi _\pm \| ^2 \leq \eta ^2, $ i.e.
   \begin{equation}
   \label{A.3} \| P \psi _\pm \| \geq \sqrt{1 - \eta ^2} > 1 - \eta,
   \end{equation}
and
   $ \big\arrowvert \langle P \psi _+ , P \psi _- \rangle - \langle \psi _+ ,
                     \psi _- \rangle \big\arrowvert \leq  
                      2 \eta + \eta ^2$, implying that
   \begin{equation}
   \label{A.4} \big\arrowvert \langle P \psi _+, P\psi _- \rangle \big\arrowvert \leq
               \theta + 2 \eta + \eta ^2 .
   \end{equation}
In order to assure that $P\psi _+$ and $P\psi _-$ are linearly independent we request 
that $ \big\arrowvert \langle P \psi _+, P \psi _- \rangle \big\arrowvert < \| P \psi _+ \|
      \| P \psi _-\|. $
In view of \eqref{A.3} and \eqref{A.4} this latter inequality is satisfied when 
$0 < \eta < \frac{1-\theta }{4}.$
But by the definition of $\eta$ and $D$, one has
   $2C \eta ^{-1} = D > \frac{8C}{1-\theta }. $ Thus we proved that 
  $P\psi _+$ and $P\psi _-$ are linearly independent, contradicting
our assumption. Hence there are at least two (counted
with multiplicities) eigenvalues of $A$ inside the circle of radius $D$ and center
$\lambda _{i_1}$.
\hspace*{\fill}$\square $


\section{Quasimodes for the bulk of ${\rm spec}(T^{\alpha , \beta}_N)$}
\label{quasibulk}

We want to apply Proposition~\ref{Lemma 4.1} (ii) to the bulk of the spectrum of $T^{\alpha ,
\beta }_N$, 
   \[ \{\lambda ^N_{2\ell - 1}, \lambda ^N_{2\ell } \,| \,\, M < \ell < N - M\}
   \]
where $M \equiv M_N = [F(N)]$. For $M < \ell < N - M$ and $N \geq 3$
arbitrary choose $\mu \equiv \mu_\ell ^N$ to be the $\ell$'th double eigenvalue of
$T^{0,0}_N$, $\mu ^N_\ell :=-2 \cos \frac{\ell \pi }{N}$.
Our construction of quasimodes follows the standard procedure of perturbation theory
of double eigenvalues:
first we construct two approximate eigenvectors $\psi ^{\ell}_{0,\pm} $ of the operator
$\prod _\ell \circ (T^{\alpha , \beta }_N - T^{0,0}_N) \big\arrowvert _{E_\ell }$ where
$E_\ell$ denotes the two dimensional eigenspace of the eigenvalue $-2 \cos \frac{\ell \pi }
{N}$ of  $T^{0,0}_N$ and the operator is the composition of the restriction of the perturbation $T
^{\alpha , \beta }_N - T^{0,0}_N$ to $E_\ell $ with the orthogonal projection
$\prod _\ell $ onto $E_\ell $. The two quasimodes $\psi ^\ell _\pm $ are then
obtained by adding a first order correction to $\psi ^{\ell}_{0,\pm}  $.
To this aim introduce
$\widetilde\psi ^\ell_+=\psi^{N,N+\ell }$ and $\widetilde\psi ^\ell_-=\psi^{N,N-\ell}$ where
we recall that $\psi ^{N,k}$ denotes the eigenvector of $T^{0,0}_N$ with eigenvalue $2 \cos \frac{k
\pi }{N}$ defined by \eqref{psiN}. One has
   \[ \widetilde\psi ^\ell_\pm = (2N)^{-1/2} \sum ^{2N-1}_{n = 0}
                     e^{\mp in\pi \frac{N-\ell }{N}} \theta ^N_n
   \]
   \[ T^{0,0}_N \widetilde\psi ^\ell _\pm = - 2 \cos \frac{\ell \pi }{N} \cdot
      \widetilde \psi ^\ell _\pm \ \mbox { and } \ \langle \widetilde \psi ^\ell _+ ,
      \widetilde \psi ^\ell _- \rangle = 0 .
   \]
Denote by $\hat\alpha_k,\hat\beta_k$, $k \in {\mathbb Z}$
the Fourier coefficients of $\alpha,\beta$ and set
\[\hat\gamma_\ell:=\hat\beta_\ell -2 \cos \frac{\ell \pi }{N}\hat\alpha _\ell,
\quad \quad e^{- i \eta _\ell }:= \hat \gamma _\ell / |\hat \gamma _\ell |
      \ \mbox { if } \ \hat \gamma _\ell \not= 0,
  \quad \mbox{and} \quad
e^{i\eta _\ell }:= 1  \mbox{ if } \hat \gamma _\ell = 0.
\]
 For any $M < \ell < N-M,$ let $\psi^\ell_\pm:=\psi ^{\ell}_{0,\pm} +\varphi^\ell_\pm$ where
   \[ \psi ^{\ell}_{0,\pm} :=\frac{\widetilde\psi^\ell_+\pm e^{i\eta_\ell}\widetilde\psi^{\ell}_-}
      {\sqrt 2} \quad \mbox {and} \quad
	\varphi^\ell_\pm := - \sum_ {n \ne N\pm \ell} \frac{\langle
      \psi^{N,n},(T^{\alpha,\beta}_N-T^{0,0}_N)\psi ^{\ell}_{0,\pm} \rangle}{2
      \cos \frac{\ell \pi }{N} + 2 \cos \frac{n \pi }{N}} \psi^{N,n}.
   \]

\begin{lemma}
\label{Lemma 5.1} The elements $\psi ^\ell_\pm,\ M < \ell < N-M$, of ${\mathcal H}_{2N}$
satisfy

\begin{description}
\setlength{\labelwidth}{9mm}
\setlength{\labelsep}{1.4mm}
\setlength{\itemindent}{0mm}

\item[{\rm (i)}] $\langle \psi ^\ell_+, \psi ^\ell_- \rangle = O(\frac {K^2_{\alpha , \beta }}{M^2})$ \quad and \quad 
$\|\psi ^\ell_\pm \| = 1+O( \frac {K^2_{\alpha , \beta }} {M^2})$;

\item[{\rm (ii)}] $\big\| \big(T^{\alpha , \beta }_N + 2 \cos \frac{\ell \pi }{N}
\big) \psi ^\ell_\pm \big\| = O \left( K^2_{\alpha , \beta } \frac {1}{N^2 M} \right)
\quad \mbox{where} \quad  K_{\alpha , \beta } = \| \alpha \| _{C^2} + \| \beta \| _{C^2} + 1.$
\end{description}
\end{lemma}

\medskip

First we need to establish the following auxiliary result.
\begin{lemma}
\label{element}
(i) For any $M < \ell < N-M$ and $n\neq N \pm \ell $ 
    \begin{equation*}
    \vert\langle \psi^{N,n},(T^{\alpha,\beta}_N-T^{0,0}_N) \psi ^{\ell}_{0,\pm} \rangle\vert
     = O \left( \frac{ K^{ }_{\alpha , \beta }}{N^2}\left(
     \mbox{min}_{\pm} \frac {1} { (N -n  \pm \ell )^2} + \frac 1 N\right) \right).
 \end{equation*}
(ii) For any $M < \ell < N-M,\,\,\, $ 
$\langle \psi ^{\ell}_{0,+},(T^{\alpha,\beta}_N-T^{0,0}_N) \psi ^{\ell}_{0,-} \rangle
     = O \left( \frac{ K_{\alpha , \beta }}{N^3} \right)$ and 
\[ \langle \psi ^{\ell}_{0,+},(T^{\alpha,\beta}_N-T^{0,0}_N) \psi ^{\ell}_{0,+} \rangle
     =  \frac{e^{-2\pi \ell^2/N}}{2N^{2}}  \Re \hat\gamma_\ell +
 O \left( \frac{ K_{\alpha , \beta }}{N^3} \right),
\]
\[
\langle \psi ^{\ell}_{0,-},(T^{\alpha,\beta}_N-T^{0,0}_N) \psi ^{\ell}_{0,-} \rangle
     = - \frac{e^{-2\pi \ell^2/N}}{2N^{2}}  \Re \hat\gamma_\ell + 
       O \left( \frac{ K_{\alpha , \beta }}{N^3} \right).
\]
\end{lemma}

{\it Proof:} By \eqref{psiN},
$
\psi ^{\ell}_{0,\pm} =\frac{\psi^{N,N+\ell}_{}\pm e^{i\eta_l}\psi^{N,N-l}_{}}{\sqrt 2} .
$
Recall that by \eqref{4}, one has for $f \equiv 1$ the identity $\psi ^{N,k} = \psi ^{N,k}_1$
and by Proposition \ref{Proposition 3.5}, for $f$ arbitrary,
$
T^{\alpha,\beta}_N\psi^{N,k}_f=\psi^{N,k}_{D^{\alpha , \beta }_k
 f}+O(\frac {K_{\alpha,\beta}} {N^3})
$
where $D^{\alpha , \beta }_k$ is given by
   \[ 2 \cos \big( k\pi /N - i(2N)^{-1} \partial _x \big)
      + (2N)^{-2}
      \big( \beta _2(x) + 2\alpha _2(x) \cos \big(k \pi /N - i(2N)^{-1}
      \partial _x \big) \big).
   \]
For $f \equiv 1$ one has $ D^{\alpha , \beta }_k 1 = 2 \cos \big( k \pi /N \big)+g(x)$ where
 $ g(x) =  (2N)^{-2} \big( \beta _2(x) + 2\alpha _2(x) \cos \big( k \pi /N\big) \big) $
and
   \[
   (T^{\alpha,\beta}_N-T^{0,0}_N)\psi^{N,k}=(T^{\alpha,\beta}_N-T^{0,0}_N)
   \psi^{N,k}_{1}=\psi^{N,k}_{g}+O\left(\frac {K_{\alpha,\beta}} {N^3}\right).
   \]
By Lemma \ref{Lemma 3.3}(ii) we have
   \[
    \langle \psi^{N,n}_1,\psi^{N,k}_g\rangle=  (2N)^{-2}
   \left(\widehat{(\beta_2)}_{k-n}+2\cos \big( k \pi /N\big) \widehat{(\alpha_2)}_{k-n} \right) 
e^{-\pi (n - k )^2/2N}.
   \]
Choosing $k = N \pm \ell,$ item (i) then follows as by assumption $\alpha, \beta \in C^2_0(\mathbb T, \mathbb R).$ To prove item (ii) note that for $n = N \pm \ell, $ one has 
$N \pm \ell - k  \in \{0, \pm 2\ell\}. $ It then follows from the definition of $e^{ i \eta _\ell}$ that
$\langle \psi^{\ell}_{0,+},\psi^{N,N+\ell}_g -  e^{ i \eta _\ell} \psi^{N,N-\ell}_g  \rangle = 0$ and
\[\langle \psi^{\ell}_{0,+},\psi^{N,N+\ell}_g +  e^{ i \eta _\ell} \psi^{N,N-\ell}_g \rangle = 
  e^{ i \eta _\ell} \left(\widehat{(\beta_2)}_{2\ell} - 2\cos \big( \ell \pi /N\big) \widehat{(\alpha_2)}_{2\ell} \right)
\frac{e^{-2\pi \ell^2/N}}{4N^2}
\]
\[
+e^{- i \eta _\ell} \left(\widehat{(\beta_2)}_{-2\ell} - 2\cos \big( \ell \pi /N\big)      
      \widehat{(\alpha_2)}_{-2\ell} \right)
     \frac{e^{-2\pi\ell^2/N}}{4N^2}
 = \frac{e^{-2\pi \ell^2/N}}{4N^{2}}  2\Re \hat\gamma_\ell . 
\]
From these and similar computations the claimed estimates follow.
\hspace*{\fill}{$\square $}

\medskip

{\it Proof of Lemma~\ref{Lemma 5.1}:} (i) First note that for any $M < \ell < N-M$ and $0 \le n \le N$
with $n \ne N \pm \ell$,
$| 2 \cos \frac{\ell \pi }{N} + 2 \cos \frac{n\pi }{N} | \ge \frac{2M \pi }{N^2}$.
Indeed, in the case where $M < \ell \le N/2$ and $0 \le k:= N - n < \ell$ one has
\[
| 2 \cos \frac{\ell \pi }{N} - 2 \cos \frac{k\pi }{N}| 
= 2\int_{\frac{k \pi }{N}}^{\frac{\ell \pi }{N}} \sin (x) dx
\ge 2\int_{\frac{k \pi }{N}}^{\frac{\ell \pi }{N}} \frac{2x}{\pi } dx
\]
 leading to the claimed lower bound. 
All other cases are treated in a similar way.
By Lemma~\ref{element} (i) one then concludes that
   \begin{equation}\label{estiphi}
   \Vert\varphi^\ell_\pm\Vert = O ( K_{\alpha,\beta}\frac 1{M} ).
   \end{equation}
On the other hand, $\psi ^{\ell}_{0,+}$ and $\psi ^{\ell}_{0,-}$ are orthogonal to each other, and both are orthogonal to $\varphi^\ell_\pm$. Hence $\langle \psi ^\ell_+, \psi ^\ell_- \rangle =
\langle \varphi ^\ell_+, \varphi ^\ell_- \rangle.$ Combined with the above estimate one gets
$\langle \psi ^\ell_+, \psi ^\ell_- \rangle = O ( \frac {K^2_{\alpha,\beta}}{M^2} ).$ Using in addition
that $\|\psi ^\ell_{0,\pm} \| = 1$ one then has 
$\|\psi ^\ell_\pm \| = 1+O( \frac {K^2_{\alpha , \beta }} {M^2})$.
(ii) We apply standard perturbation theory and write
\begin{equation}\label{4s}
T^{\alpha , \beta }_N\psi ^{\ell}_{\pm} =(T^{0 , 0}_N+(T^{\alpha , \beta }_N-T^{0 , 0}_N))
(\psi ^{\ell}_{0,\pm} +\varphi^\ell_\pm)
\end{equation}
and split the right hand side of \eqref{4s} into four parts
$$T^{0 , 0}_N \psi ^{\ell}_{0,\pm} ,\quad T^{0 , 0}_N\varphi^l_\pm,\quad (T^{\alpha , \beta }
_N-T^{0 , 0}_N)) \psi ^{\ell}_{0,\pm} ,\quad
(T^{\alpha , \beta }_N-T^{0 , 0}_N))\varphi^\ell_\pm.
$$
Note that
$
T^{0 , 0}_N \psi ^{\ell}_{0,\pm}  =-2 \cos \frac{\ell \pi }{N}\psi^{\ell}_{0,\pm}
$
and
\begin{equation*}
T^{0 , 0}_N\varphi^l_\pm=
- \sum_{n\neq N \pm \ell }2 \cos \frac{n \pi }{N}\cdot
\frac{\langle \psi^{N,n},(T^{\alpha,\beta}_N-T^{0,0}_N) \psi ^{\ell}_{0,\pm}  \rangle}
{2 \cos \frac{\ell \pi }{N} + 2 \cos \frac{n \pi }{N}}\psi^{N,n}
\end{equation*}
\[
= - \sum_{n\neq N \pm \ell}\left(1-\frac{2 \cos \frac{\ell \pi }{N}}{2 \cos
\frac{\ell \pi }{N} + 2 \cos \frac{n \pi }{N}}\right)
\langle \psi^{N,n},(T^{\alpha,\beta}_N-T^{0,0}_N) \psi ^{\ell}_{0,\pm} \rangle
\psi^{N,n}.
\]
In view of the definition of $\varphi^\ell_\pm$, this yields the identity
\[T^{0 , 0}_N\varphi^\ell_\pm=  - 2 \cos \frac{\ell \pi }{N}\varphi^\ell_\pm - 
\sum_{n\neq N \pm \ell} \langle \psi^{N,n},
(T^{\alpha,\beta}_N-T^{0,0}_N) \psi ^{\ell}_{0,\pm} \rangle \psi^{N,n}.
\]
Combined with
\[(T^{\alpha , \beta }_N-T^{0 , 0}_N) \psi ^{\ell}_{0,\pm} =\sum_{n\neq N \pm \ell}
\langle \psi^{N,n},(T^{\alpha,\beta}_N-T^{0,0}_N) \psi ^{\ell}_{0,\pm} \rangle \psi^{N,n}
+ \sum_{s \in \{ + , -\}}
\langle \psi ^{\ell}_{0,s},(T^{\alpha,\beta}_N-T^{0,0}_N) \psi ^{\ell}_{0,\pm} \rangle \psi ^{\ell}_{0,s}
\]
one gets
\[(T^{\alpha , \beta }_N-T^{0 , 0}_N) \psi ^{\ell}_{0,\pm} + T^{0 , 0}_N\varphi^\ell_\pm
= -2 \cos \frac{\ell \pi }{N}\varphi^{\ell}_{\pm} + \sum_{s \in \{ + , -\}}
\langle \psi ^{\ell}_{0,s},(T^{\alpha,\beta}_N-T^{0,0}_N) \psi ^{\ell}_{0,\pm} \rangle \psi ^{\ell}_{0,s}.
\]
By Lemma \ref{element} (ii) it then follows that
\[(T^{\alpha , \beta }_N-T^{0 , 0}_N) \psi ^{\ell}_{0,\pm} + T^{0 , 0}_N\varphi^\ell_\pm
=  -2 \cos \frac{\ell \pi }{N}\varphi^{\ell}_{\pm} +
O \left( \frac{K_{\alpha,\beta}}{N^2}\left(\frac 1 {M^2}+\frac 1 N\right) \right).
\]
Finally, the expression
\[
(T^{\alpha , \beta }_N-T^{0 , 0}_N)\varphi^\ell_\pm=
- \sum_ {n \ne N\pm \ell} \frac{\langle
      \psi^{N,n},(T^{\alpha,\beta}_N-T^{0,0}_N)\psi ^{\ell}_{0,\pm} \rangle}{2
      \cos \frac{\ell \pi }{N} + 2 \cos \frac{n \pi }{N}} (T^{\alpha , \beta }_N-T^{0 , 0}_N)\psi^{N,n}
\]
can be estimated by Lemma \ref{element} (i) to get for some constant $C \ge 1$,
\[
\| (T^{\alpha , \beta }_N-T^{0 , 0}_N)\varphi^\ell_\pm \|
\le C \sum_ {n \ne N\pm \ell}  \frac{N^2}{M}
 \frac{ K_{\alpha , \beta }}{N^2}\left(
     \mbox{min}_{\pm} \frac {1} { (N -n  \pm \ell )^2} + \frac 1 N\right) 
 \| (T^{\alpha , \beta }_N-T^{0 , 0}_N)\psi^{N,n} \|.
\]
Inspecting the proof of Lemma \ref{element} (i) one sees that 
$\| (T^{\alpha , \beta }_N-T^{0 , 0}_N)\psi^{N,n} \| =  O( \frac {K_{\alpha,\beta}} {N^2})$, 
yielding
$\| (T^{\alpha , \beta }_N-T^{0 , 0}_N)\varphi^\ell_\pm \| = O(\frac {K^2_{\alpha,\beta}} {N^2 M}).$
Combining all the above estimates, item (ii) follows.
\hspace*{\fill}{$\square $}

\bigskip

Lemma~\ref{Lemma 5.1} allows to apply Proposition~\ref{Lemma 4.1} and leads to the following
result.
\begin{proposition}
\label{Corollary 5.2} For any $N \geq 3$ and $M < \ell < N - M$,
there exists a pair of eigenvalues $\tau ^{N, \ell }_- \le  \tau ^{N,\ell }_+$ of
$T^{\alpha , \beta }_N$ satisfying
   \[ \big\arrowvert \tau ^{N,\ell }_\pm + 2 \cos \frac{\ell \pi }{N} \big\arrowvert
      = O \left(  K_{\alpha , \beta } \frac 1 {N^2M} \right)  \quad \mbox{where} 
         \quad K_{\alpha , \beta } = \| \alpha \| _{C^2} + \| \beta \| _{C^2} + 1.
   \]
For $N$ sufficiently large these pairs are separated from each other,
   \[ \cdots < \tau ^{N,\ell }_- \leq \tau ^{N,\ell }_+ < \tau ^{N,\ell + 1}_- \leq \tau
      ^{N,\ell + 1}_+ < \cdots
   \]
\end{proposition}

\smallskip

{\it Proof:} According to Lemma~\ref{Lemma 5.1}, for any $N \geq 3$ and any $M <
\ell < N - M$
   \[ \big\| \big( T^{\alpha , \beta }_N + 2 \cos \frac{\ell \pi }{N} \big) \psi ^N_\pm
      \big\| =  O\left( \frac {K^2_{\alpha , \beta }}{N^2M} \right).
\] 
By Proposition~\ref{Lemma 4.1} (ii), there are two eigenvalues $\tau ^{N,\ell }_- \leq \tau
^{N,\ell }_+$ of $T^{\alpha , \beta }_N$ satisfying
   \[ \big\arrowvert \tau ^{N,\ell }_\pm + 2 \cos \frac{\ell \pi }{N} \big\arrowvert
       = O \left( \frac{K^2_{\alpha ,\beta } }{N^2M} \right)  .
   \]
In case $\tau ^{N,\ell }_+ = \tau ^{N,\ell }_-$, the eigenvalue has multiplicity at least
two. To see that for $N$ sufficiently large, one has $\tau ^{N,\ell }_+ <  \tau ^{N,\ell + 1 }_- $, 
recall from the proof of Lemma~\ref{Lemma 5.1} that
\[
| 2 \cos \frac{\ell \pi }{N} - 2 \cos \frac{(\ell + 1)\pi }{N} | \, \ge \, \frac{2M\pi}{N^2} 
\qquad M < \ell <  N - M .
\]
Hence by choosing $N_0$ sufficiently large, the pairs of eigenvalues
$\tau ^{N,\ell }_\pm $ with $N \geq N_0$ satisfy
$ \tau ^{N,M +1}_- \leq \tau ^{N,M + 1}_+ < \tau ^{N,M + 2}_- \leq
      \tau ^{N,M + 2}_+ < \cdots < \tau ^{N,N-1-M}_- \leq \tau ^{N,N-1-M}_+ . $
\hspace*{\fill}{$\square $}


\section{Quasimodes for the edges of ${\rm spec}(T^{\alpha , \beta }_N)$}
\label{quasiedges}

In this section we want to apply Proposition~\ref{Lemma 4.1} to the two edges of the
spectrum of $T^{\alpha , \beta }_N$. They are treated in the same way, so we concentrate on
the left edge only,
    \[ \lambda ^N_0 < \lambda ^N_1 \leq \lambda ^N_2 < \ldots < \lambda ^N_{2M-1} \leq
       \lambda ^N_{2M}
    \]
where again $M \equiv M_N = [F(N)]$. For $0 \leq j \leq 2M$, choose as approximate
eigenvalue
   \begin{equation}
   \label{21bis} \mu ^{N,j}_- = - 2 + \frac{1}{4N^2} \lambda ^-_j
   \end{equation}
where $\lambda ^-_0 < \lambda ^-_1 \leq \lambda ^-_2 < \ldots \ $ are the periodic eigenvalues
of $H_- = - d^2/dx^2 + q_-$, considered on the interval $[0,1]$. Here
   \[ q_- = \beta _2 - 2\alpha _2 \quad \mbox { and } \quad 
       \alpha _2(x) = \alpha (2x), \ \beta _2(x) = \beta (2x).
   \]
 Furthermore choose as quasimodes
   \be\label{nouveau} \varphi ^{N,j}_-(z):= \psi ^{N,N}_{g^-_j}(z) = (4N)^{-1/4} \int ^1_0 g^-_j(s) \varrho
      _N(z, \frac{1}{2} + is) e^{-2N\pi s^2} ds
   \ee
where $(g^-_j)_{j \geq 0}$ is an orthonormal basis of eigenfunctions of $H_-$. 
First we need to establish bounds for $g^-_j$ and its derivatives. By the counting
lemma (cf e.g. \cite{KP}), for any $N$ with
   \begin{equation}
   \label{22} M_N > 2(1 + \| q_-\| _0) e^{\| q_-\| _0}
   \end{equation}
it follows that for $0 \leq j \leq 2M$,
   \begin{equation}
   \label{23} |\lambda ^-_j| \leq 4 \pi ^2(M + \frac{1}{2})^2 \leq 8 \pi ^2 F(N)^2 .
   \end{equation}
Recall that $K_{\alpha , \beta } = \| \alpha \| _{C^2} + \| \beta \| _{C^2} + 1$
for any $\alpha , \beta \in C^2 ({\mathbb T})$.
\begin{lemma}
\label{Lemma 6.1} For any $N$ satisfying \eqref{22} and any $\alpha , \beta \in C^2
({\mathbb T})$,

\begin{description}
\setlength{\labelwidth}{9mm}
\setlength{\labelsep}{1.4mm}
\setlength{\itemindent}{0mm}

\item[{\rm (i)}] $\| (g^-_j)'\| _0 \leq (2 K_{\alpha , \beta } + 8 \pi ^2 F(N)^2)
^{1/2};$\quad  {\rm (ii)}  $\| (g^-_j)''\| _0 \leq 2 K_{\alpha , \beta } + 8 \pi ^2 F(N)^2$;

\item[{\rm (iii)}] $\| (g^-_j)'''\|_0 \leq (2K_{\alpha , \beta } + 8 \pi ^2 F(N)^2)
^{3/2} + 2K_{\alpha , \beta }$;

\item[{\rm (iv)}] $\| (g^-_j)^{\tiny{ IV}}\| _0 \leq 3 (2K_{\alpha , \beta } + 8 \pi ^2 F(N)^2)
^2 + 2K_{\alpha , \beta } \leq 4(2K_{\alpha , \beta } + 8\pi ^2 F(N)^2)^2 .$
\end{description}
\end{lemma}
{\it Proof:} (i) Taking the inner product of
   $ - (g^-_j)'' + q_- g^-_j = \lambda ^-_j g^-_j$
with $g^-_j$, integrating by parts and using \eqref{23} and $\| g^-_j\| _0 = 1$ yields the bound (i) for $\| (g^-_j)'\| _0.$
(ii) Using again $ (g^-_j)'' = q_- g^-_j - \lambda ^-_j g^-_j$ one gets
   \[ \| (g^-_j)''\| _0 \leq \| q_-\| _0 \| g^-_j\| _0 + |\lambda ^-_j| \| g^-_j\| _0 \leq
      2K_{\alpha \beta } + 8 \pi ^2 F(N)^2 .
   \]
(iii) is obtained by deriving $(g^-_j)'' = q_- g^-_j - \lambda _j g^-_j$
and using (i). (iv) is obtained by arguing in the same way.
\hspace*{\fill}{$\square $}

\medskip

We also need bounds for $\| g^-_j\| _{C^0}$ and $\| g^-_j\| _{C^2}$. It is convenient to
formulate the result in a general form. For a real valued potential $q \in L^2({\mathbb T})$,
denote by $(f_j)_{j \geq 0}$ an orthonormal basis of periodic eigenfunctions of $H = -
d^2/dx^2 + q$ on $[0,1]$.
\begin{lemma}
\label{Lemma 6.2} {\rm (i)} The expression $\sup _{j \geq 0} \| f_j\| _{C^0}$ is bounded
uniformly on bounded sets of potentials in $L^2({\mathbb T})$.

{\rm (ii)} For any $N$ with $M = [F(N)] > 2(1 + \| q\| _0) e^{\| q\| _0}$ and any $0 \leq
j \leq 2M$
   \[ \| f''_j\| _{C^0} \leq \big( \| q\| _0 + 8\pi ^2 F(N)^2 \big) \| f_j\| _{C^0}
      \mbox { and } \| f'_j\| _{C^0} \leq 2 \big( \| q\| _{C^0} + 8 \pi ^2 F(N)^2
      \big) .
   \]
\end{lemma}
{\it Proof:} (i) It is well known that $f_0$ doesn't vanish on
$[0,1]$. As for any $j \geq 1$, $f_j$ is orthogonal to $f_0$, it has to vanish at least
once. Hence there exists $0 \leq x_j < 1$ so that $f_j(x_j) = 0$. As a consequence,
the translate $T_{x_j} f_j = f_j(\cdot + x_j)$ is a Dirichlet eigenfunction for the
translated potential $T_{x_j} q$. Note that $\| T_{x_j} q\| _0 = \| q \| _0$ and
$\| T_{x_j} f_j\| _{C^0} = \| f_j\| _{C^0}$. Therefore $\sup _{j \geq 1} \| f_j\|
_{C^0}$ is bounded uniformly on bounded sets of potentials in $L^2({\mathbb T})$
by the corresponding result for the Dirichlet problem -- see e.g. \cite[p. 35]{PT}.
It remains to bound $\| f_0\| _{C^0}$. As $\lambda _0(q)$ is never a Dirichlet eigenvalue, one has
   \[ f_0(x) = \frac{1}{c_0} \big( y_1(x,\lambda _0) + \frac{1-y_1(1,\lambda _0)}{y_2(1,
      \lambda _0)} y_2(x,\lambda _0) \big)\mbox{ where}
   \]
   \[ c_0 = c_0(q) = \big( \int ^1_0 \big( y_1(x, \lambda _0) + \frac{1 - y_1(1,
      \lambda _0)}{y_2(1,\lambda _0)} y_2 (x, \lambda _0) \big)^2 dx \big)^{1/2}
   \]
and $y_1, y_2$ are the fundamental solutions of $-y'' + qy = \lambda y$. By
\cite[p. 7]{PT}
   \[ |y_i(x, \lambda _0(q), q) | \leq e^{\| q\| _0} \quad \forall \quad 0 \leq x \leq 1
      \mbox { and } i = 1,2 .
   \]
Further, by \cite[p 18]{PT}, $y_i(x, \lambda , q)$ is a compact function of $q \in
L^2({\mathbb T})$, uniformly on bounded subsets of $[0,1] \times {\mathbb C}$. By
\cite[p 199]{KP}, $L^2({\mathbb T}) \rightarrow {\mathbb R}, q \mapsto \lambda _0(q)$
is a compact function as well and so is $q \mapsto y_2(1, \lambda _0(q), q)$. As
$\lambda _0(q)$ is never a Dirichlet eigenvalue $y_2(1, \lambda _0(q),q)
> 0$ for any $q$ in $L^2({\mathbb T})$. By the compactness, $y_2(1, \lambda _0(q),
q)$ is uniformly bounded away from $0$ on bounded sets of potentials in
$L^2({\mathbb T})$. Similarly, one argues by compactness to conclude that $c_0(q) > 0$
is uniformly bounded away from $0$ on bounded sets of potentials in $L^2({\mathbb T})$.

(ii) Note that
   $ \| f''_j\| _{C^0} \leq \big( \| q\| _{C^0} + |\lambda _j| \big) \| f_j\| _{C^0} .
   $
Hence the claimed estimate of $\| f''_j\| _{C^0}$ follows from item (i) and \eqref{22}
- \eqref{23}. Finally, for any $0 \leq x,y \leq 1$,
   $ f'_j(x) = f'_j(y) + \int ^x_y f''_j(s)ds .
   $
Integrate in $y$ and apply the Cauchy-Schwarz inequality to conclude that
   \begin{align*} \| f'_j\| _{C^0} &\leq \| f'_j\| _0 + \| f''_j\| _0 \leq \big( \| q\| _{C^0} + 8 \pi ^2 F(N)^2 \big) ^{1/2} +
                     \big( \| q \| _{C^0} + 8 \pi ^2 F(N)^2 \big)
   \end{align*}
where the latter inequality follows from the proof of Lemma~\ref{Lemma 6.1} (i), (ii).
\hspace*{\fill }$\square $

\begin{lemma}
\label{Lemma 6.3} For any $N$ with $M =[F(N)] > 2(1 + \| q_-\| _0)e^{\| q_-\|_0}$ and any
$0 \leq j \leq 2M$, the elements $\varphi ^{N,j}_-$ in ${\mathcal H}_{2N}$ satisfy
the following estimates: 
\begin{description}
\setlength{\labelwidth}{9mm}
\setlength{\labelsep}{1.4mm}
\setlength{\itemindent}{0mm}

\item[{\rm (i)}] $|\langle \varphi ^{N,j}_-, \varphi ^{N,k}_- \rangle - \delta _{j,k}|
\leq \frac{1}{4\pi N} (2K_{\alpha , \beta } + 8\pi ^2 F(N)^2) \quad \forall 0 \leq k \leq 2M.$

\item[{\rm (ii)}] $\| (T^{\alpha ,\beta }_N + 2 - \frac{1}{4N^2} \lambda ^-_j)\varphi
^{N,j}_-\| \leq \frac{F(N)^2}{N^3} (K_{\alpha , \beta } + 1)^2 C$
\end{description}
where $C > 0$ can be chosen uniformly on $L^2$-bounded subsets of $C^2({\mathbb T})$.
\end{lemma}

{\it Proof:} (i) By the definition \eqref{4} and Lemma~\ref{Lemma 3.3} (iii),
   $ |\langle \varphi ^{N,j}_-, \varphi ^{N,k}_- \rangle - \langle g^-_j, g^-_k \rangle |
      \leq \frac{1}{4\pi N} \| (g^-_j)'\| _0 \| (g^-_k)'\| _0 .
   $
By Lemma~\ref{Lemma 6.1} (i), we get $\| (g^-_j)'\| _0 \| (g^-_k)'\| _0 \leq 2 K_{\alpha , \beta }
+ 8 \pi ^2 F_1(N)^2$ and hence the claimed estimate.
(ii) By the triangle inequality
   \begin{align}
   \begin{split}
   \label{30} &\| T^{\alpha , \beta }_N \varphi ^{N,j}_- - \mu ^{N,j}_- \varphi ^{N,j}_-\|\leq \| T^{\alpha , \beta }_N \varphi ^{N,j}_- - \psi ^{N,N}_{D^{\alpha ,
                 \beta }_N(g^-_j)} \| + \| \psi ^{N,N}_{D^{\alpha , \beta }_N(g^-_j)}
                 - \mu ^{N,j}_- \varphi ^{N,j}_- \| .
   \end{split}
   \end{align}
Let us begin by estimating the latter term. By definition, $\varphi ^{N,j}_- = \psi ^{N,
N}_{g^-_j}$ and hence
   $\mu ^{N,j}_- \varphi ^{N,j}_- = \psi ^{N,N}_{\mu ^{N,j}_- g^-_j}.
   $
By Lemma \ref{Lemma 3.3} (iv) we then conclude
   \begin{equation}
   \label{32} \| \psi ^{N,N}_{D^{\alpha , \beta }_N(g^-_j)} - \mu ^{N,j}_-
                 \varphi ^{N,j}_-\| \leq
              \| D^{\alpha , \beta }_N(g^-_j) - \mu ^{N,j}_- g^-_j\| _0 .
   \end{equation}
As $ D^{\alpha , \beta }_\ell = 2 \cos \big( \frac{\ell\pi}{N} - \frac{i}{2N} \partial _x
                     \big) + \frac{1}{4N^2} \big( \beta _2 + \alpha _2 2 \cos (\frac{\ell \pi}{N}
                     - \frac{i}{2N} \partial _x) \big)$ one gets for $\ell = N$
$ D^{\alpha , \beta }_N - 2 \cos \big( - \frac{i}{2N} \partial _x
                     \big) + \frac{1}{4N^2} \big( \beta _2 - \alpha _2 2 \cos ( -
                     \frac{i}{2N} \partial _x) \big)$. Furthermore,
   $\mu ^{N,j}_- g^-_j = \big( - 2 - \frac{1}{4N^2} \partial ^2_x \big) g^-_j +
      \frac{1}{4N^2} \big( \beta _2 - 2 \alpha _2 \big) g^-_j.$ Hence we get
   \begin{align}
   \begin{split}
   \label{34} \| D^{\alpha , \beta }_N (g^-_j) - &(\mu ^{N,j}_-) g^-_j \| _0 \leq 
	\| 2 \cos \big( - \frac{i}{2N} \partial _x \big) g^-_j - \big(2 + 
          \frac{1}{4N^2} \partial ^2_x \big) g^-_j \| _0 \\
              &\quad + \frac{1}{4N^2} \| 2\alpha _2 \big( 1 - \cos (- \frac{i}
                 {2N} \partial _x) \big) g^-_j \| _0 .
   \end{split}
   \end{align}
The latter two terms are estimated individually. 
   \[ 2 \cos \big( - \frac{i}{2N} \partial _x \big) g^-_j = \sum _{n \in {\mathbb Z}}
      2 \cos \big( \frac{\pi n}{N} \big) \big( \widehat{g^-_j} \big) _n e^{i 2\pi nx} .
   \]
Using the Taylor expansion of $2 \cos \frac{\pi n}{N}$,
one concludes that
   \begin{align*} &\| 2 \cos \big( - \frac{i}{2N} \partial _x \big) g^-_j - \big(
                     2 + \frac{1}{4N^2} \partial ^2_x \big) g^-_j \| _0 \\
                  &\leq \frac{1}{12} \frac{1}{(2N)^4} \big( \sum _{n \in {\mathbb Z}}
                     \big\arrowvert (\widehat{g^-_j})_n\big\arrowvert ^2 (2\pi n)^8
                     \big) ^{1/2} \leq \frac{1}{12} \frac{1}{(2N)^4} \| (g^-_j)^{\tiny IV}\| _0 .
   \end{align*}
By Lemma~\ref{Lemma 6.1}, it then follows that
   \begin{align}
   \begin{split}
   \label{36} &\big\| 2 \cos \big( - \frac{i}{2N} \partial _x \big) g^-_j - \big( 2 +
                 \frac{1}{4N^2} \partial ^2_x \big) g^-_j \big\| _0 \leq \frac{1}{48} \frac{1}{N^4} \big(2 K_{\alpha , \beta } + 8\pi ^2 F(N)^2
                 \big)^2 .
   \end{split}
   \end{align}
In a similar way one estimates
   \begin{align*} &\frac{1}{4N^2} \big\| 2 \alpha _2 \big( 1 - \cos (- \frac{i}{2N} \partial _x )
                     g^-_j\big) \big\|_0 \leq \frac{\| \alpha \| _{C^0}}{4N^2} \big( \sum _{n \in {\mathbb Z}}
                     \big( \frac{\pi n}{N} \big) ^4 \big\arrowvert (\widehat{g^-_j})_n
                     \big\arrowvert ^2 \big) ^{1/2} \\
                  &\leq \frac{\| \alpha \| _{C^0}}{(2N)^4} \| (g^-_j)''\| _0 \leq \frac{K_{\alpha , \beta }}{16N^4} \big( 2K_{\alpha , \beta } +
                     8\pi ^2 F(N)^2 \big) \leq \frac{1}{32N^4} \big( 2K_{\alpha , \beta }
                     + 8 \pi ^2 F(N)^2 \big) ^2
   \end{align*}
where for the latter inequality we again used Lemma~\ref{Lemma 6.1}. Combining \eqref{32} -
\eqref{36} yields
   \begin{equation}
   \label{40} \big\| \psi ^{N,N}_{D^{\alpha , \beta }_N(g^-_j)} - \mu ^{N,j}_- \varphi ^{N,j}
              _- \big\| \leq \frac{F(N)^4}{16 N^4} (8\pi ^2 + 2K_{\alpha , \beta })^2 .
   \end{equation}
It remains to estimate the first term on the right hand side of \eqref{30}. By
Theorem~\ref{Proposition 3.5} (with $\ell = N$),
   \begin{equation}
   \label{42} \big\| T^{\alpha , \beta }_N \varphi ^{N,j}_- - \psi ^{N,N} _{D^{\alpha ,
              \beta }_N(g^-_j)} \big\| \leq \frac{1}{N^3} K_{\alpha , \beta } \| g^-_j\|
              _{C^2} .
   \end{equation}
By Lemma~\ref{Lemma 6.2}, for any $0 \leq j \leq 2M_N$
   \begin{equation}
   \label{44} \| g^-_j \| _{C^2} \leq \big( 2K_{\alpha , \beta } + 8 \pi ^2 F (N)^2 \big)
              C
   \end{equation}
where $C \geq 1$ can be chosen uniformly on bounded subsets of $L^2({\mathbb T})$.
Combining \eqref{30}, \eqref{40}, \eqref{42}, and \eqref{44} yields
   $\big\| T^{\alpha , \beta }_N \varphi ^{N,j}_- - \mu ^{N,j}_- \varphi
                     ^{N,j}_- \big\| \leq \frac{F(N)^4}{16N^4} (8 \pi ^2 + 2K_{\alpha , \beta })^2 +
                     \frac{F(N)^2}{N^3} K_{\alpha , \beta } (2K_{\alpha , \beta } + 8\pi
                     ^2) C$
where $C \geq 1$ can be chosen uniformly on bounded subsets of $L^2({\mathbb T})$.
\hspace*{\fill }$\square $

\begin{proposition}
\label{Corollary 6.4} For any $N \geq 3$ and any $0 \leq j \leq 2M$ there exists an
eigenvalue $\tau ^{N,j}_-$ of $T^{\alpha , \beta }_N$ satisfying
   \[ \big\arrowvert \tau ^{N,j}_- - \mu ^{N,j}_-\big\arrowvert \leq \frac{F(N)^2}{N^3}
      (K_{\alpha , \beta } + 1)^2 C
   \]
where $C > 0$ can be chosen uniformly on $L^2$-bounded subsets of $C^2({\mathbb T})$.
For $N$ sufficiently large, the eigenvalues $\tau ^{N,j}_-$ can be listed in increasing
order
   \[ \tau ^{N,0}_- < \tau ^{N,1}_- \leq \tau ^{N,2}_- < \ldots < \tau ^{N,2M-1}_-
      \leq \tau ^{N,2M}_N .
   \]
\end{proposition}

{\it Proof:} According to Lemma~\ref{Lemma 6.3}, for any $N \geq 3$ and $0 \leq j,k \leq
2M$,
   \[ \big\| \big( T^{\alpha , \beta }_N - \mu ^{N,j}_- \big) \varphi ^{N,j}_- \big\|
      \leq C_N:= \frac{F(N)^2}{N^3} (K_{\alpha , \beta } + 1)^2 C
   \]
where $C$ can be chosen uniformly on $L^2$-bounded subsets of $C^2({\mathbb T})$.
Furthermore, for $N_0$ sufficiently large,
   $\big\arrowvert \langle \varphi ^{N,j}_- , \varphi ^{N,k}_- \rangle - \delta _{j,k}
      \big\arrowvert < \frac{1}{2} \quad \forall N \geq N_0 .
   $
We now apply Proposition~\ref{Lemma 4.1} (i) or (ii) depending on whether $\mu ^{N,j}_-$ is
sufficiently isolated or not. Note the the pair $\mu ^{N,2\ell }_-, \mu ^{N,2\ell - 1}
_-$ is separated from $\{ \mu ^{N,j}_- \big\arrowvert 0 \leq j \leq 2M \} \backslash
\{ \mu ^{N,2\ell }_-, \mu ^{N,2\ell -1}_- \} $ by at least $O(N^{-2})$, uniformly on
$L^2$-bounded sets of $\alpha $'s and $\beta$'s.
If $\mu ^{N,2\ell }_- -
\mu^{N,2\ell - 1}_- \leq 2C_N$, then
   \begin{align*} \big\| \big( T^{\alpha , \beta }_N - \mu ^{N,2\ell }_- \big)
                     \varphi ^{N,2 \ell - 1}_- \big\| &\leq C_ N + \big\arrowvert \mu
                     ^{N,2\ell }_- - \mu ^{N,2\ell -1 }_- \big\arrowvert \big\| \varphi
                     ^{N,2\ell - 1}_- \big\| \\
                  &\leq C_N + 2C_N \big( 1 + O(F(N)^2N^{-1}) \big)
   \end{align*}
Applying Proposition~\ref{Lemma 4.1} (ii) to
   $ \mu ^{N,2\ell }_-, \varphi ^{N,2\ell }_-, \varphi ^{N,2\ell - 1}_- \mbox { and }
      D_N = 2 \cdot 8 \cdot 3C_N
   $
we conclude that there are two eigenvalues $\tau ^{N,2\ell - 1}_- \leq \tau ^{N,2\ell }
_-$ of $T^{\alpha , \beta }_N$ so that
   \[ \big\arrowvert \tau ^{N,j}_- - \mu ^{N,2\ell }_-\big\arrowvert \leq D_N =
      \frac{F(N)^2}{N^3} (K_{\alpha , \beta } + 1)^2 C
   \]
for $j \in \{ 2\ell , 2\ell - 1 \} $, where $C$ can be chosen uniformly on $L^2$-bounded
subsets of $C^2({\mathbb T})$. If $\mu ^{N,2\ell }_- - \mu ^{N,2\ell - 1}_- > 2C_N$,
then apply Proposition~\ref{Lemma 4.1} (i) to conclude that for $j \in \{ 2\ell , 2\ell - 1
\} $, there exists an eigenvalue $\tau ^{N,j}_-$ of $T^{\alpha , \beta }_N$ so that
$|\tau ^{N,j}_- - \mu ^{N,j}_-| \leq C_N$. In particular, we then conclude that
   $\tau ^{N,2\ell - 1}_- < \tau ^{N,2\ell }_- .
   $
   
Recall that the pairs $\mu ^{N,2\ell }_- , \mu ^{N,2\ell - 1}_-$ are separated from
each other by $O(N^{-2})$. As $F(N) \leq N^\eta $ with $0 < \eta < 1/2$ it then
follows from the definition of $C_N$ that for $N$ sufficiently large
   $ \tau ^{N,0}_- < \tau ^{N,1}_- \leq \tau ^{N,2}_- < \ldots < \tau ^{N,2N - 1}
      _-  \leq \tau ^{N,2M}_-.$\hspace*{\fill }$\square $


\section{Asymptotics of the periodic eigenvalues}
\label{eigenvalues}

The aim of this section is to prove
Theorem~\ref{Theorem 1.1} stated in the introduction.

{\it Proof of Theorem~\ref{Theorem 1.1}} 
In view of
Proposition~\ref{Corollary 5.2}, Proposition~\ref{Corollary 6.4},
and the result corresponding to Proposition~\ref{Corollary 6.4}
for the right edge of the spectrum we have obtained three groups of
eigenvalues. At the left and right edge of ${\rm spec}(T^{\alpha , \beta }_N)$ there
are according to Proposition~\ref{Corollary 6.4},  $2M + 1$ eigenvalues, which for $N$
sufficiently large are different from each other when counted with multiplicities. In
the bulk of ${\rm spec}(T^{\alpha , \beta }_N)$, we found according to
Proposition~\ref{Corollary 5.2}, $N - M - 1$ pairs of eigenvalues of $T^{\alpha , \beta }
_N$ which for $N$ sufficiently large are again different from each other. It remains
to show that
   \[ \tau ^{N,2M} < \tau ^{N,2M + 1} \quad \mbox{and} \quad \tau ^{N,2N - 2M - 2} < \tau
      ^{N,2N - 2M - 1} .
   \]
To see it, note that by the Taylor expansion of $\cos$
and \eqref{23}, we have
   $$\mu ^{N,2M}_- = - 2 + \frac{1}{N^2} \lambda ^-_{2M} \leq - 2 + \pi ^2 \frac{(M+ 1/2)^2}{N^2}.
   $$
   
   Hence
   $- 2 \cos \frac{(M + 1) \pi }{N} - \mu ^{N,2M}_- \geq \frac{M\pi ^2}{N^2} +
      O \big( \frac{1}{N^2} \big) .
   $
   Moreover, by Proposition~\ref{Corollary 6.4},
   $ \tau ^{N,2M} - \mu ^{N,2M}_- = O \big( \frac{M^2}{N^3} \big)
   $
and 
   $ \tau ^{N,2M + 1} + 2 \cos \frac{(M + 1)\pi }{N} = O \big( \frac{1}{N^2} \big) 
   $ by Proposition~\ref{Corollary 5.2}.
Therefore, for $N$ sufficiently large,
  $\tau ^{N,2M} < \tau ^{N,2M + 1} .
   $
Similarly one shows that
   $ \tau ^{N,2N-2M-2} < \tau ^{N,2N-2M -1} .
   $
Hence the eigenvalues $(\tau ^{N,n})_{0 \leq n \leq 2N -1}$ of $T^{\alpha , \beta }_N$
are listed in increasing order (and with multiplicities) and thus coincide with
$(\lambda ^N_n)_{0 \leq n \leq 2N - 1}$. The claimed estimates now follow from
Proposition~\ref{Corollary 5.2} and Proposition~\ref{Corollary 6.4}.
\hspace*{\fill }$\square $
\medskip

To finish  this section, let us mention that for smooth potentials and with some effort, our method allows to compute the full asymptotic expansion in $\frac 1{N^2}$ of all the eigenvalues of $T^{\alpha,\beta}_N$. For eigenvalues in the bulk, such an asymptotic expansion is obtained by regular perturbation theory at any order. As the eigenvalues come in separated pairs and the eigenvalues forming such a pair might coincide, 
their asymptotics are obtained via a $2\times 2-$block diagonalization and a subsequent straightforward diagonalization of the (symmetric) $2\times 2-$blocks. For eigenvalues in one of the edges of the spectrum,
the asymptotics is obtained by adding corrections to the 'densities' $g_j^\pm$ in \eqref{nouveau},
obtained by improving Theorem \ref{Proposition 3.5} so that the remainder term can be chosen to be
of arbitrary order in $N^{-2}$. These corrections can be explicitly computed by solving homological  equations, obtained from the asymptotic expansion of the operator $D^{\alpha , \beta }_\ell $ in \eqref{nouveau2} and solved by inverting certain Hill operators. 
As the eigenvalues of a Hill operator come in separated pairs and two eigenvalues forming such a pair might coincide, one obtains their asymptotics also via a $2\times 2-$block diagonalization.


\section{Asymptotics of the discriminant}
\label{discriminant}

The principal goal of this section is to prove Theorem \ref{Theorem9.1} concerning the asymptotics of the discriminant $\Delta_N(\mu)$.
Recall (cf. \cite{HK1}, Section 2) that $\Delta_N^2(\mu)-4$ is related to the characteristic polynomial of $\jmq$ as follows
\begin{equation}
   \label{1.1} \Delta _N(\mu )^2 - 4 = q^{-2}_N \prod ^{2N-1}_{j=0} (
               \lambda ^N_j - \mu ) .
   \end{equation}
   First we derive asymptotics of $q_N=\prod_1^N(1+\frac 1{4N^2}\alpha(\frac nN))$. For later reference we derive at the same time also asymptotics for 
$p_N:= \frac{1}{2} tr (Q^{\alpha, \beta}_N)=
\frac12\sum_0^{2N-1}\lambda_n^N=\frac1{4N^2}\sum_1^N\beta(i/N)$. 
It turns out that the asymptotics of $p_N$ is better than one could expect from the asymptotics 
of the eigenvalues in Theorem \ref{Theorem 1.1}.
   
\begin{proposition}
\label{jacobbb} 
Uniformly on bounded subsets of functions $\alpha , \beta $ in $C_0^2({\mathbb T})$,
   \begin{equation}
   \label{jacobi}
   q_N=1+O(N^{-3})\ \ \ \mbox{ and }\ \ \ p_N=O(N^{-3}).
   \end{equation}
\end{proposition}

{\it Proof:} As $\int ^1_0 \beta (x) dx = 0$ it follows from \eqref{A.1bis} that
   \[ p_N = \frac{1}{2} tr (Q^{\alpha , \beta }_N) = (2N)^{-2} \sum ^N_{i = 1}
      \beta \big( \frac{i}{N} \big) =O(N^{-3}).
   \]
Similarly, one has $\sum ^N_{i=1} (2N)^{-2} \alpha (\frac{i}{N} ) = O(N^{-3})$ and
thus
   \[ q_N = \exp \left( \sum ^N_{i=1} \log \left( 1 + \frac{1}{4N^2} \alpha (\frac{i}
      {N} ) \right)\right) = \exp (O(N^{-3}))
   \]
leading to the claimed estimate $q_N = 1 + O(N^{-3})$.
\hspace*{\fill }$\square $

\medskip

In the introduction we have also introduced the discriminants $\Delta_\pm$.
For $q_\pm = 0$ one gets $\Delta (\lambda) = 2\cos
(\sqrt{\lambda}/2)$ (cf end of Appendix A) and hence, with $\pi _n = \pi n$ for
$n \geq 1$,
   \[ \Delta (\lambda)^2 - 4 = - 4 \sin ^2(\sqrt{\lambda }/2) = - \lambda
      \prod _{n \geq 1} \frac{(4n^2 \pi ^2 - \lambda )^2}
      {16 \pi ^4_n} .
   \]
Similarly (see \cite{KP}), for $q_\pm $ arbitrary, one gets the following product representation
   \begin{equation}
   \label{1.2} \Delta _\pm (\lambda )^2 - 4 = (\lambda ^\pm _0 - \lambda )
               \prod _{n \geq 1} \frac{(\lambda ^\pm _{2n} - \lambda )(\lambda
               ^\pm _{2n-1} - \lambda )}{16 \pi ^4_n} .
   \end{equation}
Finally recall from the introduction  that $\Lambda ^{\pm ,M}$ denote the boxes 
   \[ \Lambda ^{\pm ,M}  \equiv \Lambda^{\pm , M}_2:= [\lambda ^\pm _0 - 2, \lambda
      ^\pm_{2[F(M)]} + 2] + i[-2,2]
   \]
   where $M=[F(N)]$ and $N\geq N_0$. We chose
$N_0 \in {\mathbb Z}_{\geq 1}$  so that
   \be\label{81bis1} \lambda ^\pm_{2k+1} - \lambda ^\pm _{2k} \geq 6 \quad \forall \
   k \geq F(F(N_0)).
   \ee

%


The estimates \eqref{9.1} and \eqref{9.2} of Theorem \ref{Theorem9.1} are
obtained in a similar fashion so we concentrate on the proof of
\eqref{9.1} only. We first need to establish several
auxiliary results.  
We cover $\Lambda^{-,M}$ by open neighborhoods, each containing one spectral band and
its adjacent gaps,
$$
\Lambda^{-,M}=\bigcup_{n \le [F(M)]}\Lambda^-_{n,\rho}
$$
where $ \Lambda^-_{1,\rho }:= [\lambda ^-_0 - 3, \lambda ^-_2 + 2\rho ] + i[-3,3]$ and for 
$2\leq n \le [F(M)] -1$,
   \[ \Lambda^-_{n,\rho}:= [\lambda ^-_{2n-3} - 2\rho , \lambda ^-_{2n} + 2\rho ]
      + i [-3,3] \qquad \mbox{and} \qquad
   \]
   \[ \Lambda ^-_{[F(M)],\rho }:= [\lambda ^-_{2[F(M)]-3}-2\rho ,  \lambda ^-_{2[F(M)]} + 3 ] + i[-3,3]
   \]
with $\rho > 0$ chosen so that
   \begin{equation}
   \label{81bisA} \lambda ^\pm _{2k} + 2\rho < \lambda ^\pm _{2k+1} - 2\rho \quad
                 \forall \ k \geq 0 .
   \end{equation}
We will study the asymptotics  of $\Delta^2_N \big(-2+\frac{1}{4N^2}\lambda\big)$
in each domain $\Lambda^-_{n,\rho}$ separately.
For this purpose  introduce 
\[
P^N_{1}(\mu):=\varepsilon^{-3}\Pi_{j=0}^2(\lambda^N_{j}-\mu) \quad \mbox{and} \quad
 P^N_{n}(\mu):=\varepsilon^{-4}\Pi_{j=2n-3}^{2n}(\lambda^N_{j}-\mu) \,\,\, \forall 2 \leq n \leq M
 \]
 where $\varepsilon=1/4N^{-2}$. Furthermore define for  $1\leq n\leq M$, $Q^N_n(\mu)$ 
\[ \Delta^2_N(\mu)-4 = \varepsilon^4 P^N_n (\mu) Q^N_n(\mu).
  \]
Defining $\pi_k=k\pi$ for $k\neq 0$ and $\pi_0=1$,  we write similarly
$\Delta^2_-(\lambda)-4 = \frac 1 {16\pi^4} P^-_1 (\lambda) Q^-_1(\lambda)$ with
\[
   P^-_1 (\lambda)=\Pi_{0\leq j\leq 2}(\lambda^-_j-\lambda) \quad \mbox{and} \quad
 Q^-_1(\lambda)=\Pi_{k
\geq 2}\frac{(\lambda^-_{2k}-\lambda)(\lambda^-_{2k-1}-\lambda)}{16\pi_k^4}
\]
 whereas for
$2\leq n \leq M,$ we define
$\Delta^2_-(\lambda ) - 4 = \frac{1}{16\pi ^4_{n-1}}\frac{1}{4\pi ^2_n}P^-_n(\lambda )
Q^-_n(\lambda )$ with
  \[ P^-_n(\lambda) := (\lambda^-_{2n}-\lambda)(\lambda^-_{2n-1}
  -\lambda) (\lambda^-_{2n-2}
  -\lambda) (\lambda^-_{2n-3}
  -\lambda)
  \]
and
   \[ Q^-_n(\lambda ) = \frac{1}{4 \pi ^2_n} (\lambda ^-_0 - \lambda )
     \prod _{{k \not= n, n-1}} \ \frac{(\lambda ^-_{2k} - \lambda )
     (\lambda ^-_{2k-1} -
     \lambda )}{16 \pi ^4_k} .
   \]

By Theorem~\ref{Theorem 1.1}, for $\lambda $ in $\Lambda ^-_{n,\varrho }$ with $2 \leq n
\leq M,$

      \begin{align}\label{Da.101}
         P^N_n (-2+\frac{1}{4N^2}\lambda ) &=\prod^{2n}_{j=2n-3}\big(\frac{
         \lambda^N_j-(-2+\varepsilon\lambda^-_j)}
         {\varepsilon}+
         (\lambda^-_j-\lambda)\big)\nonumber \\
	    &=\prod^{2n}_{j=2n-3}\big(\lambda^-_j - \lambda + O\big(\frac{M^2}{N}\big)
         \big) =P^-_n(\lambda) + O\big(n^3\frac{M^2}{N}\big) 
	  \end{align}
where we used that $\lambda^-_j-\lambda=O(n)$ for $\lambda\in\Lambda^-_{n,\rho}$ and
$2n-3\leq j\leq2n$.
Similarly , for $n=1$, one has
\[
P^N_1 (-2+\frac{1}{4N^2}\lambda )= P^-_1(\lambda )+O(\frac{M^2}N).
\]
In order to prove \eqref{9.1} we show that, for $2\leq n\leq M$,
\begin{equation}
\label{9.4}
\Delta ^2_N(-2+\varepsilon\lambda)-4= \frac{1}{16\pi ^4_{n-1}}\frac{1}{4\pi ^2_n}
\big( P^-_n(\lambda)+O\big(n^3\frac{M^2}
{N}\big) \big) Q^-_n(\lambda) + O\big(\frac{M^2}{N}\big) + O\big( \frac{1}{M} \big)
\end{equation}
uniformly for $\lambda$ in $\Lambda ^-_{n,\varrho}$ and, for $n=1$,
\begin{equation}
\label{9.4bis}
\Delta ^2_N(-2+\varepsilon\lambda)-4= \frac{1}{16\pi ^4}\big( P^-_1(\lambda)+O\big(\frac{M^2}
{N}\big) \big) Q^-_1(\lambda) + O\big( \frac{M^2}{N}\big) + O\big( \frac{1}{M} \big)
\end{equation}
uniformly for $\lambda$ in $\Lambda ^-_{1,\varrho}$. The estimates
\eqref{9.4} and \eqref{9.4bis} are proven in two steps.

\begin{lemma}
\label{Lemma9.2} Uniformly for any $\mu=-2+\varepsilon\lambda$ where $\lambda \in
\Lambda ^-_{n,\varrho}$ and $1 \leq n \leq F(M)$
   \[\prod^{2N-2M-2}_{j=2M+1} (\lambda^N_j-\mu)=\frac{N^{4M+2}}{(2\pi )^{4M}(M!)^4} \big(1+
      O\big(\frac{n^2}{M}\big) \big) .
 \]
\end{lemma}

{\it Proof:} Set $\xi^N_{2N-1}=4$ and, for $1\leq \ell \leq N-1$,

\be\label{8.4ter}\xi^N_{2\ell}=\xi^N_{2\ell-1}=2(1-\cos\frac{\ell\pi}{N}).
 \ee

Note that $\xi^N_j$ is an increasing sequence satisfying for $2M+1\leq j \leq 2N-2M-2$

\[\xi^N_j\geq \xi^N_{2M+1}>\xi^N_{2M}=2\big(1-\cos\frac{M\pi}{N}\big)>\frac{M^2\pi^2}{N^2}
\big(1-\frac{\pi^2}{12}\frac{M^2}{N^2}\big).
 \]

For $j \in \{ 2\ell, 2\ell - 1 \}$ and
$\mu=-2+\varepsilon\lambda$ with $\lambda \in \Lambda ^-_{n,\varrho}$, in view of
Theorem~\ref{Theorem 1.1},
   \begin{align*} \lambda^N_j-\mu &=-2\cos\frac{\ell \pi }{N}+O
                     \big(\frac{1}{F(N) N^2}\big)+2-\varepsilon\lambda \\
                   &=\xi^N_j - \varepsilon \lambda +O\big(\frac{1}{F(N) N^2}\big)= \xi ^N_j + O\big(\frac{n^2}{N^2}\big)
   \end{align*}

where we used that $\lambda = O(n^2)$ for $\lambda \in \Lambda ^-_{n,\varrho }$. Hence
\[
\prod^{2N-2M-2}_{j=2M+1}\big(\lambda_j^N-\mu\big)=\prod^{2N-2M-2}_{j=2M+1}
\xi ^N_j \big(1+\frac{1}
{\xi _j^N}O\big(\frac{n^2}{N^2}\big)\big)
\]
\begin{equation}
\label{9.5} =\prod^{2N-2M-2}_{j=2M+1}\xi^N_j\cdot \prod^{2N-2M-2}_{j=2M+1} \big(1+\frac{1}
{\xi _j^N} O \big(\frac{n^2}{N^2}\big)\big) .
\end{equation}

The latter two products are estimated separately.
As by \eqref{8.4ter}, $\frac 1 {\xi_j^N}=O(\frac{N^2}{M^2})$ for $2M < j <
2N-2M-1$ and hence $\frac 1 {\xi_j^N}O(\frac{n^2}{N^2})=O(\frac{F(M)^2}{M^2})$ for
$1\leq n\leq F(M)$
one can estimate
\[\sum^{2N-2M-2}_{j=2M+1}\log\big(1+\frac{1}{\xi _j ^N} O\big(\frac{n^2}{N^2}\big)\big) =
O\big(\frac{n^2}{N^2}\big) \sum^{N-M-1}_{\ell=M+1}\frac{1}{1-\cos\frac{\ell\pi}{N}}.
\]
Note that $\frac{1}{1-\cos x}$ is a monotonically decreasing function on $[\frac{M\pi}{N},
\pi]$. Hence
\[\frac{\pi }{N}\sum^{N-M-1}_{\ell =M+1}\frac{1}{1-\cos\frac{\ell \pi}{N}}\leq \int\limits^{\pi-\frac{M\pi}{N}}_{\frac{M\pi}{N}}\frac{dx}{1-\cos x} .
\]
Taking into account that
$1-\cos x= 2\sin^2\frac{x}{2}
$
and making the change of variable of integration $t:=\frac{x}{2}$ we get
\[\int\limits^{\pi-\frac{M\pi}{N}}_{\frac{M\pi}{N}} \frac{dx}{1-\cos x} = \int\limits^{\pi-\frac{M\pi}{N}}_{\frac{M\pi}{N}} \frac{dx}{2\sin^2\frac{x}{2}}=
\int\limits^{\frac{1}{2} \big(\pi-\frac{M\pi}{N}\big)} _{\frac{M\pi}{2N}} \frac{dt}
{\sin^2 t}
  = -\frac{\cos t}
{\sin t}\big\arrowvert ^{\frac{1}{2}(\pi - \frac{M\pi}{N})}_{\frac{M\pi}{2N}} \leq
\frac{N}{M}.
\]
Thus
$ \sum^{2N-2M-2}_{j=2M+1}\log\big(1+\frac{1}{\xi _j^N} O\big(\frac{n^2}{N^2}\big)\big) = O
\big(\frac{n^2}{M}\big)
$
leading to
\[\prod^{2N-2M-2}_{j=2M+1}\big(1+\frac{1}{\xi_j^N} O\big(\frac{n^2}{N^2}\big)\big) = e^{O
\big(\frac{n^2}{M}
\big)} = 1+ O\big(\frac{n^2}{M}\big).
\]
It then follows that
\[\prod^{2N-2M-2}_{j=2M+1} (\lambda^N_j-\mu) = \big(1+O\big(\frac{n^2}{M}\big)\big)
\cdot \prod^{2N-2M-2}_{j=2M+1} \xi^N_j \qquad \qquad \qquad
\]
\[
=2^{2(N-2M-1)}\prod^{N-M-1}_{\ell=M+1}\big(1-\cos\big(\frac{\ell\pi}{N}\big)
         \big)^2\cdot
         \big(1+O\big(\frac{n^2}{M}\big)\big) \qquad
\]
\[= 2^{2(N-2M-1)}\big(\frac{\prod^{N-1}_{\ell=1}\big(1-\cos\frac{\ell\pi}{N}\big)}
         {\prod^M_{\ell=1}
         \big(1-\cos\frac{\ell\pi}{N}\big)\prod^M_{\ell=1}\big(1-\cos\frac{(N-\ell)\pi}{N}
         \big)}\big)^2
         \big(1+O\big(\frac{n^2}{M}\big)\big) .
\]
By Lemma~\ref{LemmaA.1}, Lemma~\ref{LemmaA.3}, and Lemma~\ref{LemmaA.4} this latter
expression can be estimated by
\[2^{2N-4M-2}\frac{(2N2^{-N})^2\big(1+O\big(\frac{1}{N}\big)\big)}{\big(\frac{\pi^2}{2N^2}
\big)^{2M}(M!)^4
\big(1-O\big(\frac{M^3}{N^2}\big)\big)
2^{2M}\big(1-O\big(\frac{M^3}{N^2}\big)\big)}\big(1+O\big(\frac{n^2}{M}\big)\big)
\]
\[=\frac{1}{2^{4M}}\frac{N^{4M+2}}{\pi^{4M}(M!)^4}\big(1+O\big(\frac{n^2}{M}\big)\big)
\]
which is the claimed estimate.
\hspace*{\fill }$\square$

\begin{lemma}
\label{Lemma9.3} Uniformly for $\mu=-2+\varepsilon\lambda$ with $\lambda$ in $\Lambda ^-
_{n,\varrho }$
and $1 \leq n \leq M$,
\[\prod^{2N-1}_{j=2N-2M-1}(\lambda^N_j-\mu)=2^{4M+2}\big(1+O\big(\frac{M^3}{N^2}\big)\big).
\]
\end{lemma}

{\it Proof:} In view of Theorem~\ref{Theorem 1.1} (right edge), for any $0\leq j\leq 2M,\ \lambda \in\Lambda ^-_{n,\varrho }$,
\[
\lambda^N_{2N-1-j}-\mu =2-\varepsilon\lambda^+_j + O \big(\frac{M^2}{N^3}
                \big)-(-2+\varepsilon\lambda)  = 4 + O\big( \frac{M^2}{N^2} \big)
\]
as $\lambda ^+_j = O(M^2)$ for any $0 \leq j \leq 2M$. Thus
\begin{align*}
\prod^{2N-1}_{j=2N-2M-1}(\lambda^N_j-\mu) & =\prod^{2M}_{j=0}(\lambda^N_{2N-1-j}-\mu)  =4^{2M+1}\prod^{2M}_{j=0}\big(1+ O\big(\frac{M^2}{N^2}\big)\big) \\
&=2^{4M+2}\mbox{exp} \big(\sum^{2M}_{j=0}\log\big(1+O\big(\frac{M^2}{N^2}\big)\big) \big).
\end{align*}
Using the bound
$\sum^{2M}_{j=0} \log\big( 1+ O\big(\frac{M^2}{N^2}\big)\big)=O\big(\frac{M^3}{N^2}\big)
$
we get 
\[
\prod^{2N-1}_{j=2N-2M-1}(\lambda^N_j-\mu) =2^{4M+2} \mbox{exp} \big(O\big(\frac{M^3}{N^2}
\big)\big)=2^{4M+2} \big(1+O\big(\frac{M^3}{N^2}\big)\big)
\]
which is the claimed estimate.\hspace*{\fill }$\square$
\begin{lemma}
\label{Lemma9.4} Uniformly for $\mu=-2+\varepsilon\lambda$ with 
$\lambda \in \Lambda ^-_{n, \varrho }$, the product
$\prod^{2M}_{j=0}(\lambda^N_j-\mu)$ satisfies the following estimates:
(i) for $2\leq n\leq F(M)$
\[\frac{\pi^{4M}}{4} \frac{(M!)^4}{N^{4M+2}}\frac{1}{16 \pi ^4_{n-1}} \frac{1}{4 \pi ^2_n}\big( P^-_n(\lambda)+O
(\frac{n^3M^2}{N})
\big)\cdot Q^-_n(\lambda) \cdot (1+O(\frac{n^2}{M})) \cdot (1 + O(\frac{M^2}{N}));
 \]

(ii) for $n=1$
\[\frac{\pi^{4M}}{4} \frac{(M!)^4}{N^{4M+2}}\frac{1}{16 \pi ^4} \big( P^-_1(\lambda)+O
(\frac{M^2}{N})
\big)\cdot Q^-_1(\lambda) \cdot (1+O(\frac{1}{M})) \cdot (1 + O(\frac{M^2}{N})).
 \]
\end{lemma}

{\it Proof:} In view of the asymptotics of Theorem~\ref{Theorem 1.1}, with $\varepsilon=
(2N)^{-2}$,
\begin{align*}
\lambda^N_j-\mu & =-2+\varepsilon\lambda^-_j+2-\varepsilon\lambda+O(\frac{M^2}{N^3})=\varepsilon(\lambda^-_j-\lambda+O(\frac{M^2}{N})) .
\end{align*}
Hence
$ \prod^{2M}_{j=0}(\lambda^N_j-\mu)  =\varepsilon^{2M+1}\prod^{2M}_{j=0}(\lambda^-_j
-\lambda +O(\frac{M^2}{N})).
$
The items $(i)$ and $(ii)$ are proved in a very similar way - in fact $(ii)$ is a
little simpler. Hence we concentrate on $(i)$, i.e. the case where $2\leq n\leq F(M)$.
Given $\lambda $ in $\Lambda ^-_{n,\varrho }$, the latter product is
split up into three parts,
\[\prod^{2n-4}_{j=0}(\lambda^-_j-\lambda+O(\frac{M^2}{N})) \cdot \prod^{2n}_{j=2n-3}(\lambda^-_j
- \lambda+ O(\frac{M^2}{N})) \cdot
\prod^{2M}_{j=2n+1}(\lambda^-_j-\lambda+O(\frac{M^2}{N})).
\]
Note that
\[\prod^{2n}_{j=2n-3}(\lambda^-_j-\lambda+O(\frac{M^2}{N}))=P^-_n(\lambda)+O(\frac{n^3M^2}{N})
\]
uniformly for $\lambda \in \Lambda ^-_{n,\varrho}$ with $2 \leq n \leq M$, cf. equation \eqref{Da.101}. Next
consider
\[\prod^{2n-4}_{j=0} \big( \lambda^-_j -\lambda+ O (\frac{M^2}{N}) \big) =\prod^{2n-4}_{j=0}(\lambda^-_j-\lambda)  \cdot \prod^{2n-4}_{j=0} \big( 1+\frac{1}{\lambda^-_j - \lambda}O(\frac{M^2}{N}) \big).
\]
We claim that
\[\prod^{2n-4}_{j=0}(1+\frac{1}{\lambda^-_j-\lambda}O(\frac{M^2}{N}))=1+O(\frac{M^2}{N})
\]
uniformly for $\lambda$ in $\Lambda^-_{n,\rho}$. Indeed, by the choice of $\varrho $, the factors
$(\lambda - \lambda _{j }^-)^{-1},\ 2\ell -1\leq j\leq 2\ell,$  can be estimated
by $O((n^2 - \ell ^2)^{-1})$ uniformly for $\lambda \in \Lambda ^-_{n,\varrho }$ and $1 \leq \ell
\leq n - 2$.
Hence $\sum\limits^{2n-4}_{j=0}\frac{1}{\lambda-\lambda^-_j}$ is a bounded analytic function of $\lambda\in\Lambda^-_{n,\rho}$ with a bound depending on $\rho $.
Similarly one treats
\[\prod^{2M}_{j=2n+1}(\lambda^-_j-\lambda + O(\frac{M^2}{N}))=\prod^{2M}_{j=2n+1}(\lambda^-_j-
\lambda) \cdot \prod^{2M}_{j=2n+1}
(1+\frac{1}{\lambda^-_j -\lambda}O(\frac{M^2}{N})).
\]
Again, by the choice of $\varrho $ and the asymptotics of $\lambda ^-_j$,
$ \sum ^{2M}_{j = 2n + 1} \frac{1}{\lambda _j^- - \lambda }$ is a bounded analytic function of
$\lambda \in \Lambda^- _{n,\rho }$. Thus
   \[ \prod^{2M}_{j=2n+1} \big( 1 + \frac{1}{\lambda _j^- - \lambda } O(\frac{M^2}{N}) \big) =
      1 + O (\frac{M^2}
      {N}) .
   \]
Combining the estimates obtained we have
   \begin{align*} &\prod ^{2M}_{j = 0} (\lambda ^N_j - \mu ) = \varepsilon ^{2M + 1} \cdot
       \prod ^M_{k = 1} 2^4 \pi ^4_k \cdot (\lambda ^-_0 - \lambda ) \cdot
	   \prod ^{n-2}_{k=1} \frac{(\lambda ^-_{2k} - \lambda )(\lambda ^-_{2k-1} - \lambda )}
       {2^4 \pi ^4_k}\cdot \frac{1}{16 \pi ^4_{n-1}} \\
       &\cdot\frac{1}{16 \pi ^4_{n}}\cdot  \big( P^-_n(\lambda ) + O(\frac{n^3M^2}{N}) \big)
       \cdot \prod ^M_{k = n + 1} \frac{(\lambda ^-_{2k} - \lambda )
       (\lambda ^-_{2k-1} - \lambda )}{2^4 \pi ^4_k} \
       \cdot (1 + O(\frac{M^2}{N})) .
   \end{align*}
Finally we note that, with $\lambda ^-_j = 4k^2 \pi ^2 + \alpha _j$ for $j \in \{ 2k, 2k - 1\}$,
where $\alpha _j = O(1)$,
   \begin{align}\nonumber
   \begin{split}
   &\prod ^\infty _{k = M + 1} \frac{ (\lambda ^-_{2k} - \lambda )(\lambda ^-
                  _{2k-1} - \lambda )}
                  {2^4 \pi ^4_k} = \mbox{exp } \left\{ \sum ^\infty _{k = M + 1}
                  \log \big( \frac{\lambda ^-_{2k}-\lambda}{4 \pi ^2_k} \big) + \log
                  \big( \frac{\lambda ^-_{2k-1} - \lambda }{4\pi ^2_k} \big) \right\} \\
                  &= \mbox{exp } \left\{ \sum ^\infty _{k = M + 1} \mbox{log } \big( 1 +
                  \frac{\alpha _{2k}-\lambda }{4 \pi ^2_k}
                  \big) + \mbox{log } \big( 1 + \frac{\alpha  _{2k-1} - \lambda }{4\pi ^2_k}
                  \big) \right\} = 1 + O(\frac{n^2}{M})
   \end{split}
   \end{align}
uniformly for $\lambda \in \Lambda^- _{n,\rho }, 1 \leq n \leq M$. Here we used
   $\sum ^\infty _{k =  M + 1} \frac{1}{k ^2} \leq \int ^\infty _M \frac{1}{x^2} dx = \frac{1} {M} .
   $
Furthermore
   $ \prod ^M_{k = 1} 2^ 4 \pi ^4_k = = (2\pi ) ^{4M}(M!)^4 .
   $
By the definition of $Q^-_n(\lambda )$, we then obtain that
$\prod ^{2M}_{j = 0} (\lambda ^N_j - \mu )$ equals
   \[ \varepsilon ^{2M+1}
      \frac{(2\pi) ^{4M} (M!)^4}{16\pi ^4_{n-1}\cdot 4\pi ^2_n}\big( P^-_n
      (\lambda ) + O(\frac{n^3M^2}{N})
	  \big) Q^-_n(\lambda ) \big( 1 + O (\frac{M^2}{N}) \big) \big( 1 + O
				  (\frac{n^2}{M}) \big)
   \]
as claimed.
\hspace*{\fill}$\square$

\begin{lemma}
\label{Lemma 8.1bis} Uniformly for $\lambda \in \Lambda ^-_{n,\varrho}$, $1 \leq n
\leq M$, $Q^-_n(\lambda )=O(n^2)$.
\end{lemma}

{\it Proof:} By the Counting Lemma (cf \cite{KP}) for periodic eigenvalues there
exits $n_0 \geq 1$ so that $|\lambda ^\pm _n - 4n^2 \pi ^2| \leq 1$ for any $n \geq
n_0$. Note that $n_0$ can be chosen uniformly for bounded sets of functions $\alpha ,
\beta \in C_0^2$.
It turns out that the cases $1 \leq n < n_0$ and $n_0 \leq n
\leq M$ have to be treated separately.
However they can be proved in a similar way and so we concentrate on the case
$n_0 \leq n\leq M$ only.

\be\label{8.7n} Q^-_n(\lambda ) = \frac{1}{4 \pi ^2_n} (\lambda ^-_0 - \lambda )
     \prod _{{k \not= n, n-1}} \ \frac{(\lambda ^-_{2k} - \lambda )
     (\lambda ^-_{2k-1} -
     \lambda )}{16 \pi ^4_k}
   \ee
and that $\frac {\sin(\sqrt\lambda/2)}{\sqrt\lambda/2}$ can be written as an infinite product,
\[
\frac {\sin(\sqrt\lambda/2)}{\sqrt\lambda/2}=\prod_{m\geq 1}
\frac {m^2\pi^2-\lambda/4}{m^2\pi^2}=\prod_{m\geq 1}\frac {4\pi_m^2-\lambda}
{4\pi^2_m}.
\]
Hence  for $\lambda\in \Lambda ^-_{n,\varrho}$,
\be\label{8.7nbis}
Q^-_n(\lambda ) = \frac{\lambda_0^--\lambda}{4\pi_n^2}\left(\frac {\sin(\sqrt\lambda/2)}{\sqrt\lambda/2}\right)^2
\left(\frac {4\pi_n^2\cdot 4\pi_{n-1}^2 }{(4\pi_n^2-\lambda)(4\pi_{n-1}^2-\lambda)}\right)^2f_n^-(\lambda)
\ee
where
$ f_n^-(\lambda)=\prod_{k\neq n, n-1}\frac{(\lambda^-_{2k}-\lambda)(\lambda^-_{2k-1}
-\lambda)}{(4\pi_k^2-\lambda)^2}.
$
Clearly, uniformly for $\lambda \in \Lambda ^-_{n,\varrho}$, $n_0 \leq n
\leq M$, one has
$
\frac{\lambda-\lambda_0^-}{4\pi_n^2}=1+O(\frac 1{n})
$
and
\[
 \left(\frac {\sin(\sqrt\lambda/2)}{\sqrt\lambda/2}\right)^2
\left(\frac {4\pi_n^2\cdot 4\pi_{n-1}^2 }{(4\pi_n^2-\lambda)(4\pi_{n-1}^2-\lambda)}\right)^2=
\qquad \qquad
\]
\[
\left(\frac {\sin(\sqrt\lambda/2)}{(\pi_n-\sqrt\lambda/2)(\pi_{n-1}-\sqrt\lambda/2)}\right)^2\cdot 4^3\cdot \frac{\pi_n^4\pi_{n-1}^4}{\lambda(2\pi_n+\sqrt\lambda)^2(2\pi_{n-1}+\sqrt\lambda)^2}=O(n^2)
\]
where we used that for $\lambda \in \Lambda ^-_{n,\varrho}$, $n_0 \leq n
\leq M$,
$
\frac{\sin{(\sqrt\lambda/2)}}{(\pi_n-\sqrt\lambda/2)(\pi_{n-1}-\sqrt\lambda/2)}=O(1).
$
Finally we need to estimate $f_n^-(\lambda)$. For $n \geq n_0$, by
the choice of $\rho>0$ there exists $\rho'>0$ so that
$
\vert 4\pi_k^2-\lambda\vert\geq\frac 1 {\rho'}\vert k^2-n^2\vert,\ \ \ \forall k\neq n,n-1,\ \ \forall\lambda\in \Lambda ^-_{n,\varrho}.
$
Thus
\begin{eqnarray}\nonumber
\left\vert\frac{(\lambda^-_{2k}-\lambda)(\lambda^-_{2k-1}-\lambda)}{(4\pi_k^2-\lambda)^2}\right\vert &\leq
\left(1+\left\vert\frac{\lambda^-_{2k}-4\pi_k^2}{4\pi_k^2-\lambda}\right\vert\right)
\left(1+\left\vert\frac{\lambda^-_{2k-1}-4\pi_k^2}{4\pi_k^2-\lambda}\right\vert\right)\qquad \\
&\leq
\left(1+\rho'\left\vert\frac{\lambda^-_{2k}-4\pi_k^2}{k^2-n^2}\right\vert\right)
\left(1+\rho'\left\vert\frac{\lambda^-_{2k-1}-4\pi_k^2}{k^2-n^2}\right\vert\right)
\end{eqnarray}

Using that $\sum_{k\geq 1}k^{-2}\leq\pi^2/6$, one has by Cauchy-Schwarz
\[
\sum_{k\neq n,n-1}\left\vert\frac{\lambda^-_{2k}-4\pi_k^2}{k^2-n^2}\right\vert,\quad
\sum_{k\neq n,n-1}\left\vert\frac{\lambda^-_{2k-1}-4\pi_k^2}{k^2-n^2}\right\vert\leq\pi K
\quad \mbox{ where }\]
\[K:=\big(\sum_{k\geq 1} \vert\lambda^-_{2k}-4\pi_k^2\vert^2
+\vert \lambda^-_{2k-1}-4\pi_k^2 \vert^2\big)^{\frac 1 2}.
\qquad \qquad \quad
\]
Hence, uniformly for $\lambda \in \Lambda ^-_{n,\varrho}$, $n_0 \leq n
\leq M$,
\[
\vert f_n^-(\lambda)\vert = \prod_{k\neq n, n-1}\frac{(\lambda^-_{2k}-\lambda)
(\lambda^-_{2k-1}-\lambda)}{(4\pi_k^2-\lambda)^2} \qquad \qquad
 \]
\[\leq  \exp\left(\sum_{k
\neq n,n-1}\log{\big(1+\rho'\left\vert\frac{\lambda^-_{2k}-4\pi_k^2}
{k^2-n^2}\right\vert\big)}+\sum_{k\neq n,n-1}\log{\big(1+\rho'\left\vert\frac{\lambda^-_{2k-1}
-4\pi_k^2}{k^2-n^2}
\right\vert\big)}\right)
 \]
\[\leq \exp (2 \rho '\pi K). \qquad \qquad \qquad \qquad
\]
Altogether, $Q^-_n(\lambda )=O(n^2)$ uniformly for
$\lambda \in \Lambda ^-_{n,\varrho}$, $n_0 \leq n \leq M$ as claimed.
\hspace*{\fill }$\square $

\bigskip

{\it Proof of Theorem~\ref{Theorem9.1}:} By Proposition~\ref{jacobbb}
the factor $q_N^{-2}$ appearing in the product representation \eqref{1.1}
of $\Delta_N^2-4$ satisfies the asymptotics
   \[ q^{-2} _N = 1 + O\big( \frac{1}{N^3} \big) .
   \]
Combining Lemma~\ref{Lemma9.2},
Lemma~\ref{Lemma9.3}, and Lemma~\ref{Lemma9.4} one obtains, uniformly for $\lambda $
in $\Lambda ^-_{n,\rho }$ with $1 \leq n \leq F(M)$
   \[ \Delta ^2_N(- 2 + \frac{1}{4N^2} \lambda ) - 4 = \big[
                  \Delta ^2_-(\lambda ) - 4 + O(\frac{M^2}{n^3N}) Q^-_n(\lambda )
				  \big] \big(1 + O ( \frac{n^2}{M}) \big) \big( 1 + O
				  (\frac{M^2}{N}) \big) .
   \]
As $Q^-_n(\lambda ) = O(n^2)$ (Lemma~\ref{Lemma 8.1bis}) and $\Delta ^2_-(\lambda )-4 = O(1)$
uniformly for $\lambda \in \Lambda ^-_{n,\varrho }$ with $1 \leq n \leq
F(M)$ it then follows that
   \[ \Delta ^2_N(-2 + \frac{1}{4N^2} \lambda ) - 4 = \Delta ^2_-(\lambda ) - 4 + O(\frac{F(M)^2}
      {M}) .
   \]
To determine how the signs of $\Delta_N$
and $\Delta_-$ are related note that for $\lambda $ in
the set
\begin{equation*}
\{z\in\mathbb{C} \mid \quad \mbox{dist}(z,
             [\lambda^-_{2n-1},\lambda^-_{2n}])< 2\varrho\},\quad  1\leq n\leq F(M),
\end{equation*}
one has
   \[ \Delta _N(-2 + \frac{1}{4N^2} \lambda ) = (-1)^{N-n} \sqrt[+]{\Delta ^2_N
      (-2 + \frac{1}{4N^2} \lambda )}
\quad \mbox{and} \quad
   \Delta _-(\lambda ) = (-1)^n \sqrt[+]{\Delta ^2_-(\lambda )}.
   \]
Hence
   \[ \Delta _N(-2 + \frac{1}{4N^2} \lambda ) = (-1)^N \Delta _-(\lambda ) +
      O(\frac{F(M)^2}{M}).
   \]
The estimates for $\lambda $ in $\Lambda^+_{n,\varrho }$ with $1 \leq n \leq F(M)$ are
obtained in a similar fashion.
Finally, to see that these estimates are uniform on bounded sets of $\alpha , \beta $
in $C_0^2({\mathbb T}, {\mathbb R})$ it suffices
to note that $\rho $ of \eqref{81bisA} can be chosen uniformly on such sets as the
periodic eigenvalues of $-\partial ^2_x + q_\pm $ are compact functions of $\alpha ,
\beta $ -- see \cite{KP}, Proposition B.11).
\hspace*{\fill }$\square $

\medskip

As $\Delta _N(\mu )$ and $\Delta _-(\lambda )$ are analytic functions one can
apply Cauchy's theorem to deduce from Theorem~\ref{Theorem9.1} corresponding
estimates of the derivatives $\partial ^j_\mu \Delta _N $ or equivalently
$\partial ^j_\lambda \Delta _N(- 2 + \frac{\lambda }{4N^2}) = \frac{1}{(4N^2)^j}
\partial ^j_\mu \Delta _N \left( -2 + \frac{\lambda }{4N^2} \right)$ as well as
$\partial ^j_\lambda \Delta _N \left( 2 - \frac{\lambda }{4N^2} \right) =
\frac{(-1)^j}{(4N^2)^j} \partial ^j_\mu \Delta _N \left( 2 - \frac {\lambda }
{4N^2} \right)$. Let
   \[ \Lambda ^{\pm , M}_1 = [\lambda^+_0 - 1, \lambda ^\pm _{2[F(M)]} + 1] +
      i [-1,1] .
   \]

\begin{corollary}
\label{Corollary 8.7} Let $F$ satisfy {\rm (F)}, $M=[F(N)]$ with $N \ge N_0$,
and  $\alpha , \beta \in C_0^2({\mathbb T}, {\mathbb R})$. Then, for any $j \geq 1$ and uniformly for $\lambda $ in $\Lambda ^{-,M}_1$,
   \[ \frac{1}{(4N^2)^j} \partial ^j_\mu \Delta _N\left( -2 + \frac{1}{4N^2} \lambda
      \right) = (-1)^N \partial ^j_\lambda \Delta _-(\lambda ) + O \left( \frac
      {F(M)^2}{M} \right)
   \]
and similarly, for $\lambda $ in $\Lambda ^{\pm , M}_1$
   \[ \frac{(-1)^j}{(4N^2)^j} \partial ^j_\mu \Delta _N\left( 2 - \frac{1}{4N^2} \lambda
      \right) = \partial ^j_\lambda \Delta _+(\lambda ) + O \left( \frac
      {F(M)^2}{M} \right) .
   \]
These estimates hold uniformly on bounded sets of functions $\alpha , \beta $ in
$C_0^2({\mathbb T}, {\mathbb R})$.
\end{corollary}

{\it Proof:} By Cauchy's theorem, for $j \geq 1$,
   \[ \frac{1}{(4N^2)^j} \partial ^j_\mu \Delta _N\left( - 2 + \frac{1}{4N^2} \lambda
      \right) = \frac{1}{j!} \frac{1}{2\pi i} \int _{\partial \Lambda ^{-,M}}
      \frac{\Delta _N(-2 + \frac{1}{4N^2} z)}{(z-\lambda )^{1 + j}} dz
   \]
and
    \[ \partial ^j_\lambda \Delta _{-}(\lambda ) = \frac{1}{j!} \frac{1}{2\pi i}
       \int _{\partial \Lambda ^{-,M}}
      \frac{\Delta _-(z)}{(z-\lambda )^{1 + j}} dz
   \]
where $\partial \Lambda ^{-,M}$ denotes the boundary of the rectangle $\Lambda
^{-,M} \equiv \Lambda^{-,M}_2$ with counterclockwise orientation. Hence
   \[ \frac{ \partial ^j_\mu \Delta _N\left( -2 +  \lambda/4N^2
      \right)}{(4N^2)^j} - (-1)^N\partial ^j_\lambda \Delta _- (\lambda ) = \frac{1}{j!} \frac{1}
      {2\pi i} \int\limits _{\partial \Lambda^{-,M}_2} \frac{\Delta _N(-2 +
      \frac{1}{4N^2}z) - (-1)^N \Delta _-(z)}{(z - \lambda )^{j+1}} dz .
   \]
For $\lambda $ in $\Lambda ^{-,M}_1, |z-\lambda |^2 \geq 1$ and hence by
Theorem~\ref{Theorem9.1}, uniformly on $\Lambda^{-,M}_1$
   \begin{align*} \frac{ \partial ^j_\mu \Delta _N\left( -2 +  \lambda/4N^2
      \right)}{(4N^2)^j} - \partial ^j_\lambda \Delta _- (\lambda )
                     &= \frac{1}{j!} \frac{1}{2\pi i} \int\limits _{\partial \Lambda^{-,M}_2}
                     \frac{\Delta _N(-2 +
                     \frac{1}{4N^2}z) - \Delta _-(z)}{(z - \lambda )^{j+1}} dz \\
                  &= O \left( \frac{F(M)^2}{M} \right) .
   \end{align*}
By this argument, also the uniformity statement with respect to $\alpha , \beta $
follows.
\hspace*{\fill }$\square $

\medskip

Corollary~\ref{Corollary 8.7} allows to obtain asymptotics of the zeroes 
 of $\dot \Delta _N(\mu ):=\frac d{d\mu} \Delta _N(\mu )$ at the edges in terms of the zeroes
of $\dot \Delta _\pm (\lambda ):=\frac d{d\lambda}\Delta _\pm (\lambda )$. One sees in a
 straightforward way that the $N-1$ zeroes of the polynomial $\dot \Delta _N(\mu )$ are all real and simple and when listed in increasing order, satisfy 
$\lambda^N_{2n-1} \le \dot \lambda ^N_n \le \lambda^N_{2n}$ 
for any $0 < n < N.$ 
Similarly one sees that the zeroes of $\dot \lambda ^\pm _n$ are all real and simple, and when listed in increasing order, satisfy $\lambda^\pm_{2n-1} \le \dot \lambda ^\pm_n \le \lambda^\pm_{2n}$ for any $n \ge 1$.

\begin{corollary}
\label{Corollary 8.8} Let $F$ satisfy {\rm (F)}, $M=[F(N)]$,
and  $\alpha , \beta \in C_0^2({\mathbb T}, {\mathbb R})$. Then for any $1 \leq n \leq F(M)$,
   \[ \dot \lambda ^N_n = - 2 + \frac{\dot \lambda ^-_n}{4N^2} + O \left(
      \frac{n^2}{N^2} \frac{F(M)^2}{M} \right)
\,\,\, \mbox{ and } \,\,\,\,
   \dot \lambda^N_{N-n} = 2 - \frac{\dot \lambda ^+_n}{4N^2} + O \left(
      \frac{n^2}{N^2} \frac{F(M)^2}{M} \right) .
   \]
These estimates hold uniformly on bounded sets of functions $\alpha , \beta $ in
$C_0^2({\mathbb T}, {\mathbb R})$.
\end{corollary}

{\it Proof:} The asymptotics of the zeroes of $\dot \Delta _N(\mu )$ at the two edges are obtained in a similar fashion so we concentrate
on the ones at the left edge. Let $\Gamma ^-_n$ be the contour of the box 
$[ \lambda^-_{2n-1} - \rho , \lambda ^-_{2n} + \rho ] +i [-1,1]$,  contained in $\Lambda^{-,M}_1$,
   \[ \Gamma ^-_n = \partial ([ \lambda^-_{2n-1} - \rho , \lambda ^-_{2n} + \rho ]
      +i [-1,1])
   \]
where $\rho$ is chosen as in (\ref{81bisA}).  By Theorem~\ref{Theorem 1.1},  for $N$ sufficiently large, $\dot \lambda _n^N$ 
is the only zero of $\dot \Delta _N(\mu )$ in the box 
$-2 + \frac{1}{4N^2}([ \lambda^-_{2n-1} - \rho , \lambda ^-_{2n} + \rho ] +i [-1,1])$.
In particular,  $\dot \Delta _N(\mu )$ doesn't vanish on the contour 
$\Gamma ^N_n = -2 + \frac{1}{4N^2}\Gamma ^-_n$.
By Cauchy's theorem it then follows that for any $1 \leq n \leq F(M)$,
 $1 = \frac{1}{2\pi i} \int _{\Gamma ^N_n} \frac{\partial ^2_\mu \Delta _N}
      {\partial _\mu \Delta _N} d \mu
  $
and $ \dot \lambda ^N_n = \frac{1}{2\pi i} \int _{\Gamma ^N_n} \mu
      \frac{\partial ^2_\mu \Delta _N}{\partial _\mu \Delta _N} d\mu.
   $
Hence
\[ \dot \lambda ^N_n = - 2 + \frac{1}{2\pi i} \int _{\Gamma ^N_n}
                     (\mu + 2) \frac{\partial^2_\mu \Delta _N(\mu)}{\partial _\mu \Delta _N(\mu)})
                     d\mu 
\]                 
and with the change of variable $\mu = -2 + \frac{\lambda }{4N^2}$
   \[ 4N^2 (\dot \lambda ^N_n + 2) =  \frac{1}{2\pi i} \int _{\Gamma ^-
      _n} \lambda \frac
		{\partial ^2_\lambda \Delta _N \left( -2 +
			\frac{\lambda } {4N^2} \right) }
		{\partial _\lambda \Delta _N \left( -2 +
			\frac{\lambda } {4N^2} \right) } d\lambda .
   \]
Similarly one has
   $ \dot \lambda ^-_n = \frac{1}{2\pi i} \int _{\Gamma ^-_n} \lambda
	\frac {\partial ^2_\lambda \Delta _-(\lambda )}
		{\partial _\lambda \Delta _- (\lambda )} d \lambda .
   $
The difference $4N^2 (\dot \lambda ^N_n + 2) - \dot \lambda ^-_n $ thus equals
   \begin{align*} &\frac{1}{2\pi i} \int _{\Gamma ^-_n} \lambda \left( \frac{\partial
                     ^2_\lambda \Delta _N \left( -2 + \frac{\lambda }{4N^2}\right)}
                     {\partial _\lambda \Delta _N \left( -2+ \frac{\lambda }{4N^2}
                     \right) } - \frac{(-1)^N \partial ^2_\lambda \Delta _-(\lambda )}
                     {(-1)^N \partial _\lambda \Delta _-(\lambda )} \right) d \lambda \\
                  &= \frac{1}{2\pi i} \int _{\Gamma ^-_{n}} \lambda \frac{ \partial ^2
                     _\lambda \Delta _N \left( -2 + \frac{\lambda }{4N^2} \right) -
                     (-1)^N \partial ^2_\lambda \Delta _-(\lambda )}{\partial_\lambda
                     \Delta _N \left( -2 + \frac{\lambda }{4N^2} \right) } d\lambda \\
                  &+ \frac{1}{2\pi i} \int _{\Gamma ^-_n} \lambda \frac{\partial ^2
                     _\lambda \Delta _-(\lambda ) \cdot \left( (-1)^N \partial _\lambda \Delta _-(\lambda )
                     - \partial _\lambda\Delta _N\left( -2 + \frac{\lambda }{4N^2} \right) \right)} 
                     {\partial _\lambda \Delta _N \left(
                     -2 + \frac{\lambda }{4N^2} \right) \partial _\lambda \Delta _-
                     (\lambda )} d \lambda.
   \end{align*}
The two latter integrals are estimated separately. Use Corollary~\ref{Corollary 8.7} and the facts that on $\Gamma ^-_n$, $\lambda = O(n^2)$ and
   \[\partial ^2 _\lambda \Delta _-(\lambda ), \quad
	\frac{1}{\partial _\lambda \Delta _N \big( -2 + \frac{\lambda }{4N^2} \big)} ,\quad
      \frac{1}{\partial _\lambda \Delta _-(\lambda ) } \quad = O(1)
   \]
to conclude that each of the two integrals is
$O \left( n^2 \frac{F(M)^2}{M} \right) $, yielding
   \[4N^2 (\dot \lambda ^N_n + 2) = \dot \lambda ^-_n   + O \left( n^2 \frac{F(M)^2}{M} \right) .
   \]
The statement on the uniformity of the estimates is obtained by using that a corresponding one for the discriminants and their derivatives holds.
\hspace*{\fill }$\square $

\medskip

\begin{appendix}
\section{Auxiliary results}
\label{Appendix A}

In this appendix we prove auxiliary results needed to compute the
asymptotics of the discriminant.

\begin{lemma}
\label{LemmaA.1} For $N \rightarrow \infty $,
   \begin{equation}
       \label{A.1} \prod ^{N-1}_{n=1} \big( 1 - \cos \frac{n\pi }{N} \big)
                  = 2N2^{-N} (1 + O(N^{-1})).
    \end{equation}
\end{lemma}

{\it Proof:} Note that $1 - \cos (n\delta ) > 0$ for $1 \leq n \leq N$,
and $\delta := \pi /N$. To compute the product in \eqref{A.1} we therefore
can take the logarithm, yielding,
   \[ \sum ^N_{n = 1} \log (1 - \cos (n\delta )) = \sum ^N_{n = 1} \log
      \big( \frac{1 - \cos (n\delta )}{(n\delta )^2} \big) + \sum ^N_{n = 1}
	  \log (n\delta )^2 .
   \]
Clearly
   \[ \sum ^N_{n = 1} \log (n\delta )^2 = N \log \delta ^2 + \log (N!)^2 = \log
      \big( (\frac{\pi }{N})^N N! \big) ^2 .
   \]
To compute the asymptotics of $\sum ^N_{n = 1} \log \big( \frac{1 - \cos (n
\delta )}{(n\delta )^2} \big)$ introduce
   \[ f(x) = \log \big( \frac{1 - \cos x}{x^2} \big) \quad 0 \leq x \leq
      \pi .
   \]
Note that for $0 \leq x \leq \pi $,
   \begin{align*} \frac{1 - \cos x}{x^2} &= \frac{1}{2} - \frac{1}{4!} x^2 +
                  \ldots = \frac{1}{2} (1 - \frac{1}{12} x^2 + \ldots ) > 0 .
   \end{align*}
Hence $f(x)$ is a well-defined, smooth function on the interval $[0,\pi ]$. Now
apply the well known formula for approximating the sum $\sum ^N_{n = 1} f(n\delta )$
by an integral (cf e.g. \cite{AE})
   \begin{equation}
   \label{A.1bis} \sum ^N_{n = 1} f(n \delta ) = \frac{1}{\delta } \int ^\pi _0 f(x)dx +
                  \frac{f(\pi ) - f(0)}{2} + O(\delta )
   \end{equation}
where the error term $O(\delta )$ is bounded by
   \[ \delta \frac{1}{12} \sup _{0 \leq x \leq \pi } |f''(x)| \cdot
      \mbox{length} ([0,\pi ]) .
   \]
Clearly
   \[ \frac{f(\pi ) - f(0)}{2} = \frac{1}{2} \big( \log (\frac{2}{\pi ^2}) -
      \log \frac{1}{2} \big) = \log \frac{2}{\pi } .
   \]
Further, $\frac{1}{\delta } \int ^\pi _0 f(x)dx = \frac{1}{\delta } \int ^\pi
_0 (\log (1 - \cos x) - 2\log x)dx$ can be explicitly computed by Lemma~\ref{LemmaA.2},
   \[ \frac{1}{\delta } \int ^\pi _0 \log (1 - \cos x) = - \frac{\pi }{\delta }
      \log 2 = \log \frac{1}{2^N}
   \]
and
   \[ - \frac{2}{\delta } \int ^\pi _0 \log x dx = - \frac{2}{\delta } (x \log x - x)
      \big\arrowvert ^\pi _0 = 2N + \log \frac{1}{\pi ^{2N}} .
   \]
Combining all these estimates yields
   \[ \sum ^N_{n = 1} \log (1 - \cos (n\delta )) = \log \big( \frac{2}{\pi } \frac{1}
      {2^N} \frac{1}{\pi ^{2N}} \big) + 2N + \log \big( \frac{\pi ^N}{N^N} N!
	  \big) ^2 + O(\frac{1}{N})
   \]
or
   \[ \prod ^N_{n = 1} \big( 1 - \cos (\frac{n\pi }{N}) \big) = \frac{2}{\pi }
      \frac{e^{2N}}{2^N N^{2N}}(N!)^2 \big( 1 + O(\frac{1}{N}) \big) .
   \]
By Stirling's formula,
   $ N! = \sqrt{2\pi N} N^N e^{-N} \big( 1 + O(\frac{1}{N}) \big)
   $
it follows
   \[ \prod ^N_{n = 1} \big( 1 - \cos \frac{n\pi }{N} \big) = \frac{4N}{2^N}
      \big( 1 + O(\frac{1}{N}) \big)
   \]
and as $\big( 1 - \cos \frac{n\pi }{N} \big) \big\arrowvert _{n=N} = 2$ we then conclude that
   \[ \prod ^{N - 1}_{n = 1} \big( 1 - \cos \frac{n\pi }{N} \big) = \frac{2
      N}{2^N} \big( 1 + O(\frac{1}{N}) \big).
	  \eqno{\square }
   \]

\begin{lemma}
\label{LemmaA.2} $\int ^\pi _0 \log (1 \pm \cos x)dx = - \pi \log 2$.
\end{lemma}

{\it Proof:} First note that by the change of variable of integration $x:= \pi
- s$,
   \[ \int ^\pi _0 \log (1 + \cos x)dx = \int ^\pi _0 \log (1 +
                  \cos(\pi - s))ds = \int ^\pi _0 \log(1 - \cos s) ds .
   \]
Hence, with $I:= \int ^\pi _0 \log(1 - \cos x)dx$, one has
   \begin{align*} 2I &= \int ^\pi _0 (\log (1 + \cos x) + \log (1 - \cos x))dx 
   = 2 \int ^{\pi /2 }_0 \log (\sin ^2 x)dx .
   \end{align*}
Using that $\sin ^2 x = \frac{1}{2} (1 - \cos 2x)$ and making the change of variable $s = 2x$, one 
gets
$I 
   = \frac{1}{2} \int ^\pi _0 \log (1 - \cos s)ds - \frac{\pi }{2} \log 2$ and the claim follows.
\hspace*{\fill } $\square $

\begin{lemma}
\label{LemmaA.3} For any $1 \leq M < N$,
   \begin{equation}
   \label{A.5} \big( \frac{\pi ^2}{2N^2} \big)^M (M!)^2 \geq \prod ^M_{n = 1}
               \big( 1 - \cos (\frac{n\pi }{N}) \big) \geq \big( \frac{\pi ^2}
			   {2N^2} \big) ^M (M!)^2 \exp \big( - O (\frac{M^3}{N^2})\big) .
   \end{equation}
\end{lemma}

{\it Proof:} As in the proof of Lemma~\ref{LemmaA.1}, consider the logarithm of the
product in \eqref{A.5}, to obtain, with $\delta := \pi / N$,
   \[ \sum ^M_{n = 1} \log (1 - \cos (n\delta )) = \sum ^M_{n = 1} \big( \log
      \frac{(n\delta )^2}{2} + \log (1 + \frac{2b_n}{(n\delta )^2}) \big)
   \]
where $b_n = 1 - \cos n \delta - \frac{(n\delta )^2}{2}$. Clearly
   \begin{equation}
   \label{A.6} \sum ^M_{n = 1} \log \frac{(n\delta )^2}{2} = \log \big( \frac{(M!)^2}
               {2^M} (\frac{\pi }{N}) ^{2M} \big) = \log \big( (M!)^2 (\frac
			   {\pi ^2}{2N^2})^M \big) .
   \end{equation}
Further not that $\frac{2b_n}{(n\delta )^2} < 0$ and
   \[ \Big\arrowvert \frac{2b_n}{(n\delta )^2} \Big\arrowvert = \Big\arrowvert -
      \frac{2}{4!} (n\delta)^2 + \frac{2}{6!}(n\delta )^4 - \Big\arrowvert \leq
	  \frac{1}{12} (n\delta )^2 \leq \frac{M^2}{N^2} .
   \]
As for $-1 < x < 0$,
   \[ 0 > \log (1 + x) = - \big( |x| + \frac{|x|^2}{2} + \frac{|x|^3}{3} + \ldots \big) \geq
      - |x| \frac{1}{1 - |x|}
   \]
it then follows that
   \[ 0 < - \log \big( 1 + \frac{2 b_n}{(n\delta )^2}\big) <  \frac{M^2}{N^2}
      \big(1 - \frac{M^2}{N^2} \big)^{-1} =  O\big( \frac{M^2}{N^2} \big) .
   \]
Summing up these estimates yields
   \[ 0 < - \sum ^M_{n = 1} \log \big( 1 + \frac{2b_n}{(n\delta )^2} \big) \leq
       M \cdot O ( \frac{M^2}{N^2}) =  O ( \frac{M^3}{N^2}) .
   \]
Combined with the estimate \eqref{A.6} one gets the claimed estimate.
 \hspace*{\fill } $\square $

\begin{lemma}
\label{LemmaA.4} For any $1 \leq M < N$
   \[ 2^M \geq \prod ^M_{n = 1} \big( 1 + \cos \frac{n\pi }{N} \big) \geq
      2^M \exp \big( - O (\frac{M^3}{N^2}) \big) .
   \]
\end{lemma}
		
{\it Proof:} Note that
   $1 + \cos \frac{n\pi }{N} = 2 - \frac{1}{2} \big( \frac{n\pi }{N} \big)^2
      + \ldots = 2\big( 1 - (\frac{n\pi }{2N})^2 + \ldots \big).
   $
Thus
   \[ \prod ^M_{n = 1} \big( 1 + \cos \frac{n\pi }{N} \big) = 2^M \prod ^M
      _{n = 1} \big( 1 - (\frac{n\pi }{2N} )^2 + \ldots \big) \leq 2^M
   \]
and
   \[ \prod ^M_{n=1} \big( 1 + \cos \frac{n\pi }{N} \big) =
                  2^M \exp \big( \sum ^M_{n = 1} \log \big( 1 - (\frac{n\pi }
				  {2N})^2 + \ldots \big) \big) \geq
				  \]\[\geq 2^M \exp \big( - (\frac{\pi }{2N})^2 \sum ^M_{n = 1}
				  n^2 \big) \geq 2^M \exp \big( - O (\frac{M^3}{N^2}) \big).\eqno{\square }\]	

\medskip

Finally we compute the spectral data for the operator $-d^2/dx^2$
when considered with periodic / antiperiodic boundary conditions on the interval
$[0,T]$. The fundamental solutions of $-d^2/dx^2$ are given by
$y_1(x,\lambda ) = \cos \sqrt{\lambda } x$ and $y_2(x,\lambda ) = \frac
{\sin \sqrt{\lambda } x}{\sqrt{\lambda }}$. Thus the periodic / antiperiodic
eigenvalues are
   \[ \lambda ^T_0 = 0; \ \lambda ^T_{2n} = \lambda ^T_{2n-1} = \big(
      \frac{n\pi }{T} \big) ^2 \quad \forall n \geq 1
   \]
and a basis of eigenfunctions is given by
   \[ f_0 = 1; \quad f_{2n}(x) = \cos \big( \frac{n\pi }{T} x \big) ; \quad f_{2n-1}(x) =
      \sin \big( \frac{n\pi }{T} x \big) .
   \]
The discriminant can be computed to be
   \[ \Delta _T(\lambda ) = y_1(T, \lambda ) + y'_2(T, \lambda ) = 2 \cos
      (\sqrt{\lambda } T)
   \]
hence
   \[ \Delta _T(\lambda )^2 - 4 = 4 \cos ^2 (\sqrt{\lambda } T) - 4 = - 4 \sin ^2
      (\sqrt{\lambda } T) .
   \]
As 
   $ \sin \sqrt{\mu } = \sqrt{\mu } \prod _{n \geq 1} \frac{n ^2 \pi ^2 - \mu }{n^2
      \pi ^2}$,
it then follows that
   \[ \Delta _T(\lambda )^2 - 4 
   = - 4 \lambda
      T^2 \prod _{n \geq 1} \big( \frac{n^2\pi ^2 - \lambda T^2}{n^2 \pi ^2}
	  \big) ^2
    = - 4T^2 \lambda \prod _{n \geq 1} \big( \frac{
      \frac{n^2 \pi ^2}{T^2} - \lambda }{\frac{n^2\pi ^2}{T^2}} \big)^2 .
   \]
In view of the values of $\lambda ^T_n$, it follows that
   \[ \Delta _T(\lambda )^2 - 4 = - 4T^2 \lambda \prod _{n \geq 1} \frac{(\lambda
      ^T_{2n} - \lambda )(\lambda ^T_{2n-1} - \lambda )}{(\frac{n\pi }{T})^4}.
   \]
Furthermore we compute the entire functions $\psi_k^T(\lambda),\  k\geq 1,\ $leading to
the normalized differentials $\frac{\psi_k^T(\lambda)}{\sqrt{\Delta_T^2(\lambda)-4}}
d\lambda$ characterized by
\[
\frac1{2\pi}\int_{\Gamma_n^T}\frac{\psi_k^T(\lambda)}{\sqrt[c]{\Delta_T^2(\lambda)-4}}d\lambda=\delta_{n,k}\ \ \ \forall n,k\geq 1
\]
where, as usual, $\Gamma_n^T$ is a counterclockwise contour around $\lambda_{2n}^T=\lambda^T_{2n-1}$, so that all other eigenvalues $\lambda_k^T,\ k\neq 2n,2n-1$, are in the exterior of $\Gamma_n^T$. We claim that
\[
\psi_n^T(\lambda)=c_k^T\prod_{l\neq k}\frac{\sigma_l^{T,k}-\lambda}{\left(\frac{l\pi}T\right)^2}
\quad \mbox{ with } \quad
\sigma_l^{T,k}=\left(\frac{l\pi}T\right)^2 \,\mbox{ and } \,c_k^T=\frac{2T^2}{k\pi}.
\]
Indeed, as $\sigma_l^{T,k}$ is in the $\ell $'th gap interval, it follows that
$\sigma_l^{T,k}=\left(\frac{l\pi}T\right)^2,\  \forall l\neq k.$
The constant $c^T_k$ is then determined by
$
1=\frac 1{2\pi}\int_{\Gamma_k^T}\frac{\psi_k^T(\lambda)}{\sqrt[c]{\Delta_T^2(\lambda)-4}}d\lambda.
$
As
\[
\frac{\psi_k^T(\lambda)}{\sqrt[c]{\Delta_T^2(\lambda)-4}}=
c_k^T\frac{1}{i2 T\sqrt\lambda}\frac{\left(\frac{k\pi}T\right)^2}{\lambda-\lambda_{2k}^T}
\]
one gets by Cauchy's Theorem that
$
c_k^T=\frac{2T^2}{k\pi}
$
as claimed.
In the special case where $T = 1/2$ one gets
   \[ \Delta (\lambda )^2 - 4 = - \lambda \prod _{n \geq 1} \frac{(\lambda
       _{2n} - \lambda )(\lambda _{2n-1} - \lambda )}{(2n\pi )^4}	
	   \]
where $\lambda _n \equiv \lambda ^T_n \Big\arrowvert _{T = 1/2}$ for any $n
\geq 0$ and
  $\Delta (\lambda ) \equiv \Delta _T(\lambda )\Big\arrowvert _{T = 1/2}$.
    For the entire functions $\psi_k(\lambda):=\psi^T_k(\lambda)|_{T=\frac 1 2}$ one gets
   \[
   \psi_k(\lambda)=\frac1{2k\pi}\prod_{l\neq k}\frac{(2l\pi)^2-\lambda}{(2l\pi)^2}
   \quad \mbox{ and }   \quad
   c_k:=c_k^T|_{T=\frac12}=\frac1{2k\pi}.
   \]
\end{appendix}


\begin{thebibliography}{99}

\bibitem{AE} H. Amann, J. Escher: {\em Analysis II}, Birkh\"auser, Basel, 1999.

\bibitem{BTP} D. Bambusi, T. Kappeler, T. Paul: {\em De Toda \`a KdV}, C. R. Acad. Sci. Paris, Ser I, 347(17-18), 1025-1030, 2009.

\bibitem{bkp2} D. Bambusi, T. Kappeler, T. Paul: {\em  Dynamics of periodic Toda chains with a large number of particles,} J. of Differential Equations 258, 4209-4274, 2015.

\bibitem{BP} D. Bambusi, A. Ponno: {\em On metastability in FPU}, Comm. Math. Phys.
264, 539-561, 2006.

\bibitem{BCP} G. Benettin, H. Christodoulidi, A. Ponno: {\em The Fermi-Pasta-Ulam problem and its underlying integrable dynamics}, J. Stat. Phys. 144, 793, 2001.

\bibitem{BGG} L. Berchialla, L. Galgani, A. Giorgilli: {\em Localisation of energy in FPU chains}
DCDS-A 11, 855 - 866, 2005.

\bibitem{BI} G. Berman, F. Izrailev: {\em The Fermi-Pasta-Ulam problem: 50 years of progress}, 
Chaos 15, 015104, 2005.

\bibitem{BKL} J. Biello, P. Kramer, Y. Lvov, {\em Stages of energy transfer in the FPU model,} Dynamical systems and differential equations (Wilmington, NC, 2002), DCDS Suppl., 113 -122, 2003. 

\bibitem{BGPU} A. Bloch, F. Golse, T. Paul, A. Uribe: {\em Dispersionless Toda and
Toeplitz operators}, Duke Math. J. 117, 157-196, 2003.

\bibitem{CGG} A. Carati, L. Galgani, A. Giorgilli: {\em The Fermi-- Pasta--Ulam problem 
as a challenge for the foundations in physics,} Chaos, 15, 015105, 2005.

\bibitem{C} C. Cercigniani: {\em Solitons. Theory and application}, Riv. Nuovo Cim. 7, 429-469, 1977.

\bibitem{FFM} W. Ferguson, H. Flaschka, D. McLaughlin: {\em Nonlinear Toda modes for the Toda chain,}
J. Comput. Phys. 45, 157-209, 1982.

\bibitem{FPU} E. Fermi, J. Pasta, S. Ulam: {\em Los Alamos report LA -- 1940}, in E. Fermi, {\em Collected Papers}, vol 2, 977-988, 1955.

\bibitem{F} H. Flaschka: {\em The Toda lattice I. Existence of integrals}, Phys. Rev.
B 9, 1924-1925, 1974.


\bibitem{FMM} E. Fucito, F. Marchesoni, E. Marinari, G. Parisi, L. Peliti, S. Ruffo, A. Vulpiani: {\em Approach to equilibrium in a chain of nonlinear oscillators},  J. de Physique 43, 707-713, 1982.

\bibitem{G} G. Gallavotti (Ed.): {\em Fermi--Pasta--Ulam problem: a status report}, Lect. Notes Phys., vol 728,
Springer, 2007. 

\bibitem{Gie1} D. Gieseker: {\em The Toda hierarchy and the KdV hierarchy}, Comm. Math. Phys. 181,
587-603, 1996.

\bibitem{Gie2} D. Gieseker: {\em Toda and KdV},  J. of Differential Geometry 64, 171-246, 2003.

\bibitem{HK1} A. Henrici, T. Kappeler: {\em Global action-angle variables for the
periodic Toda lattice}, IMRN 2008; Vol 2008: article ID run 031, 52 pages, doi:10.1093/imrn/rnn031.



\bibitem{HK4} A. Henrici, T. Kappeler: {\em Results on normal forms for FPU chains}, 
Comm. Math. Phys. 278, 14-177, 2008.

\bibitem{HK5} A. Henrici, T. Kappeler: {\em Resonant normal forms for even periodic FPU chains}, 
 J. European Math. Soc. 11, 1025-1056, 2009.


\bibitem{KP} T. Kappeler, J. P\"oschel: {\em KdV \& KAM}, Springer, Heidelberg, 2003.


\bibitem{LPR} R. Livi, M. Pettini, S. Ruffo, M. Sparpaglione, A. Vulpiani: {\em Relaxation to different stationary states in the Fermi--Pasta--Ulam model}, Phys. Rev. A 28, 3544-3552, 1983.

\bibitem{LPV} R. Livi, M. Pettini, S. Ruffo, A. Vulpiani: {\em Further results on the equipartition threshold in large nonlinear Hamiltonian systems}, Phys. Rev. A 31, 2741-2742, 1985.


\bibitem{CP} C. Morosi, L. Pizzocchero: {\em On the continuous limit of
  integrable lattices. I. The Kac-Moerbeke system and KdV theory},
  Commun. Math. Phys.  180 (2), 505–528, 1996. 

\bibitem{PU} T. Paul, A. Uribe: {\em A construction of quasimodes using coherent
states}, Ann. H. Poincar\'e Phys. Th\'eo. 59 (4), 357--382, 1993.

\bibitem{PL} M. Pettini, M. Landolfi: {\em Relaxation properties and ergodicity breaking in nonlinear Hamiltonian dynamics}, Phys. Rev. A 41, 768-783, 1990.

\bibitem{PT} J. P\"oschel, E. Trubowitz: {\em Inverse spectral theory}, Academic Press,
Boston, 1987.


\bibitem{SchW} G. Schneider, C. Wayne: {\em Counter-propagating waves on fluid surfaces
and the continuum limit of the Fermi-Pasta-Ulam model}, in International Conference on
Differential Equations, Vol. 1.2 (Berlin 1999), 390-404, World Sci. Publ. 2000.

\bibitem{T} M. Toda: {\em Theory of nonlinear lattices}, 2nd edition, Springer Series
on Solid-State Sciences 20, Springer, 1989.

 
\bibitem{ZK} N. Zabusky, M. Kruskal: {\em Interaction of solitons in
a collisonsless plasma and the recurrence of initial states}, Phys. Rev.
Lett. 15, 240 - 245, 1965.


\end{thebibliography}
\end{document}